\documentclass[a4paper,12pt]{ijamas}

\usepackage{amsmath,latexsym,amssymb,graphicx,epsf,amsfonts,amsthm}
\usepackage{epsf}
\usepackage{color}
\numberwithin{equation}{section}

\newtheorem{thm}{Theorem}
\newtheorem*{ass}{Assumption}
\newtheorem{defi}[thm]{Definition}

\begin{document}


\title{Structure Assisted NMF Methods for Separation of Degenerate
 Mixture Data with Application to NMR Spectroscopy}

\author{Yuanchang Sun
\thanks{Department of Mathematics and Statistics, Florida International University, Miami FL 33199, USA.}
\and Kai Huang \footnotemark[1]
\and Jack Xin
\thanks{Department of Mathematics, University of California at Irvine, Irvine, CA 92697, USA.}
}

\date{}
\maketitle

\begin{abstract}
In this paper, we develop structure assisted nonnegative matrix factorization
(NMF) methods for blind source separation of degenerate data.
The motivation originates from nuclear magnetic resonance (NMR) spectroscopy,
where a multiple mixture NMR spectra are recorded to identify chemical
compounds with similar structures. Consider the linear mixing model (LMM),
we aim to identify the chemical compounds involved when the mixing process
is known to be nearly singular. We first consider a class of data with
dominant interval(s) (DI) where each of source signals has dominant peaks over others.
Besides, a nearly singular mixing process produces degenerate mixtures.
The DI condition implies clustering structures in the data points.
Hence, the estimation of the mixing matrix could be achieved by data clustering.
Due to the presence of the noise and the degeneracy of the data,
a small deviation in the estimation may introduce errors
in the output.  To resolve this problem and improve robustness
of the separation, methods are developed in two aspects.
One is to find better estimation of the mixing matrix by allowing a
constrained perturbation to the clustering output,
and it can be achieved by a quadratic programming.
The other is to seek sparse source signals by exploiting the DI condition,
and it solves an $\ell_1$ optimization.
If no source information is available,
we propose to adopt the nonnegative matrix factorization approach by
incorporating the matrix structure (parallel columns of the mixing
matrix) into the cost function and develop  multiplicative iteration
rules for the numerical solutions.  We present experimental
results of NMR data to show the performance and reliability
of the method in the applications arising in NMR spectroscopy.
\end{abstract}



\section{Introduction}

Blind source separation (BSS) is a major area of research in signal and image processing.
It aims at recovering source signals from their mixtures without detailed knowledge of the mixing process.  BSS has been playing a central role in a wide range of signal and image processing problems such as speech recognition, sound unmixing, image separations, and text mining, to name a few 
\cite{Choi}, \cite{Cic}, \cite{Comon1}.
The goal of this paper is to study BSS methods for nearly degenerate mixtures arising from nuclear magnetic resonance spectroscopy (NMR).  Being one of the preeminent imaging techniques in chemistry,  NMR spectroscopy is frequently used by scientists to determine the molecular structures of organic compounds.  We shall consider multiple NMR spectra acquired from a mixture of chemical compounds.  Each compound has a unique spectral fingerprint defined by the number, intensity and locations of its NMR peaks.  However, when the compounds (component molecules) have similar functional groups, the peaks overlap in the composite
NMR spectra making it difficult to identify the compounds involved.  This makes the data analysis hopeless unless we can unmix or separate the mixed data into a list of source components (source spectra).  The simplest model is the linear mixing model (LMM)

\begin{equation}
\label{LMM}
X = A\,S\;, \mathrm{with}\, A_{ij} \geq 0\;, S_{ij}\geq 0\;,
\end{equation}
where $X\in \mathbb{R}^{m\times p}, A \in \mathbb{R}^{m\times n}, S\in \mathbb{R}^{n\times p}$.  Rows of $X$ represents the measured mixed signals, rows of $S$ are the source signals. The $X, S$ are sampled functions of an acquisition variable which may be time, frequency, position, wavenumber, etc. depending on the underlying physical process.  The objective of BSS is to solve for $A$ and $S$ given $X$.  If $P$ is a permutation matrix and $D$ an invertible diagonal matrix, one can immediately notice that $ AS = (APD)(D^{-1}P^{- 1}S)$,  hence $(A,S)$ and $(APD,D^{-1}P^{-1}S)$ are considered equivalent solutions in BSS.

There have been mainly two classes of BSS methods for solving (\ref{LMM}).
The first class of methods belong to statistical regime.
Among others, independent component analysis (ICA \cite{Cic})
is the most well studied statistical BSS approach, it decomposes
a mixed signal into additive source components based on the mutual
independence of the non-Gaussian source signals.  The statistical
independence requires uncorrelated source signals, and this condition
however is not always satisfied by realistic data.  For example,
the statistical independence does not hold in the NMR spectra of
chemical compounds where molecules responsible for each source
share common structural features.   The second class of methods are deterministic approaches
which include nonnegative matrix factorization (NMF) and geometrical methods.
Firstly introduced by Paatero and Tapper \cite{NMF0}, later popularized by Lee and Seung \cite{Lee},
NMF has become the prevalent method for solving nonnegative BSS problems.
NMF seeks a factorization of $X$ into product of two nonnegative matrices
by minimizing the cost function of a certain distance or divergence metric \cite{Choi}.
NMF does not impose source independence, however, some additional constraints such as sparsity and/or smoothness of
the sources and/or the mixing matrix, are often incorporated into NMF to control
the non-uniqueness and improve the physical meaning of the solutions.

Geometrical BSS methods are based on convex geometry of the data matrix $X$.
The columns of $X$ are non-negative linear combinations of those of $A$.
In the hyperspectral unmixing (HSI) setting, a condition called pure pixel
assumption (PPA) was proposed in
\cite{Chang_07}
which requires the presence
in the data of at least one pure pixel of each endmember (source signal).
Under this assumption, several approaches such as vertex component analysis (VCA)
\cite{VCA},
pixel purity index (PPI)
\cite{Boardman_93}, and N-findr
\cite{Winter_99}
have exploited geometric features of hyperspectral mixtures to determine
the smallest simplex (convex cone) containing the data.
In NMR spectroscopy, PPA was reformulated by Naanaa and Nuzillard
\cite{NN05}.
The source signals are only required to be non-overlapping at
some locations of acquisition variable (NNA).
This condition was applied to NMR data unmixing and
led to a major success of a convex cone method.  Such a local sparseness
condition has greatly reduce this problem to a convex one which is solvable by linear programming. Another geometric method is to find a simplex (convex cone) of minimum volume enclosing the data set \cite{MVT,MVSA}.
This method amounts to solving a non-convex minimization problem by finding a matrix with minimum volume under a constraint.

The different spectra come from Fourier transform of NMR measurement of absorbance of radio frequency radiation by receptive nuclear spins of the same mixture sample at different time segments when exposed to high magnetic fields.  The NMR spectra are nonnegative.  Besides, NMR spectra of different chemical compounds are usually not independent, especially as compounds (component molecules) have similar functional groups, the peaks overlap in the composite NMR spectra making it difficult to identify the compounds involved.   The classical method of Naanaa and Nuzillard (NN) requires the condition that NMR source signals to be non-overlapping at certain locations while they are allowed to overlap with each other elsewhere.  The NN assumption (NNA) requires the source signals to be strictly non-overlapping at some locations of acquisition variable (e.g., frequency).  In other words, each source signal must have a stand-alone peak where other sources are strictly zero there.   Such a strict sparseness condition leads to a dramatic mathematical simplification of a general nonnegative matrix factorization problem (\ref{LMM}) which is non-convex.  Geometrically speaking, the problem of finding the mixing matrix $A$ reduces to the identification of a minimal cone containing the columns of mixture matrix $X$.  The latter can be achieved by linear programming.  In this paper, we consider how to separate the data if NN condition is not satisfied, and the mixing process is known to be degenerate.  We are concerned with the regime where source signals do not have stand-alone peaks yet one source signal dominates others over certain intervals of acquisition variable.  In other words, a dominant interval(s) condition (DI) is required for source signals.  This is a reasonable condition for many NMR spectra.  For example, the DI condition holds well in the NMR data which motivated us.  The data is produced by the so-called DOSY (diffusion ordered spectroscopy) experiment where a physical sample of mixed chemical compounds in solvent (water) is prepared.  DOSY tries to distinguish the chemicals based on variation in their diffusion rates.   However, DOSY fails to separate them if the compounds have similar chemical functional groups (i.e., they have similar diffusion rates).  In this application, the diffusion rates of the chemicals serve as the mixing coefficients.  This presents an additional mathematical challenge due to the near singularity of the mixing matrix.  Separating these degenerate data is intractable to NN method.  New methods need to be invented.  Examination the DI condition reveals a great deal about the geometry of the mixtures.  Actually, the scattered plot of columns of $X$ must contain several clusters of points, and these clusters are centered at columns of $A$.   Hence, the problem of finding $A$ boils down to the identification of the clusters and their centroids, and it can be accomplished by data clustering methods such as K-means.  K-means clustering is computationally fast and easy to implement.  Although the data clustering in general produces a fairly good estimate of the mixing matrix, its output deviates from the true solution due to the presence of the noise, initial guess of the clustering algorithm, and so on.   In the case of nearly singular mixing matrix, a small perturbation can lead to considerable errors in the source recovery (e.g., negative spurious peaks).  To overcome this difficulty and improve robustness of the separation, we propose two remedies.  One is to find better estimation of mixing matrix by allowing a constrained perturbation to the clustering outputs, and it is achieved by a quadratic programming.  The intention is to drive the estimates closer to the true solution.  The other is to seek sparse source signals by exploiting the DI condition.  An $\ell_1$ optimization problem is formulated for recovering the source signals.

As a non-convex optimization method, NMF suffers mainly from being sensitive to
the initial guess and non-unique solutions.  These drawbacks however do not prevent
NMF from gaining its popularity in application such as spectroscopy, chemometrics,
remote sensing, image processing, and environmental science due to its capability
of reduced representation of the original data and its easy implementation.
With the limited information (other than the nonnegativity) of the data,
NMF remains as the first choice among the separation methods.
The past several decades have seen considerable efforts in the development of variants of NMF,
for example, sparse NMF, orthogonal NMF, local NMF, etc., \cite{Choi} to steer the NMF towards
the desired solutions.  One major disadvantage with NMF is the loss of
the data geometric structure in the recovery. For instance, the closeness
of the columns of $X$ should be preserved or revealed in the factorization
which a general NMF usually fails to achieve.
Some remedy has been proposed for the preservation of such geometric relationship.
In \cite{ManiNMF} an affinity graph is introduced to steer the NMF towards
the implicit geometrical information and seek a matrix factorization which
respects the graph structure.  In many NMR data cases, it is very likely that only the degeneracy of the mixing matrix is known, but no knowledge of the source signals is available. Then the NMF type methods should be used to seek a solution. In the case of degenerate mixing matrices, we will introduce additional constraints to NMF cost function to steer the solution towards the ground truth.  This case scenario that the mixing matrix $A$ is
known to be nearly degenerate (similar columns) can also be found in
in computer vision where the two adjacent frames of a moving scene
should not vary much (hence similar).

The paper is outlined as follows;  In section 2, we shall review the essentials of NN approach, then we state the new dominant interval condition on the source signals motivated by NMR data, we also discuss two cases of degenerate mixing matrix.  In section 3, we present the methods, the first method is based on the convex geometry and cluster structure the data; while the second approach belongs to NMF method where additional constraints are incorporated into the cost function to control the degeneracy of the mixing matrix.  In section 4, we
illustrate our method with numerical examples including the processing of an experimental DOSY NMR data set.
Concluding remarks are in section 5.  We shall use the following notations throughout the paper.  The notation $A^j$ stands for the $j$-th column of matrix $A$, $S^j$ for the $j$-th column of matrix $S$, $X^j$ the $j$-th column of matrix $X$.  While $S_j$ and $X_j$ are the $j$-th rows of matrix $S$ and $X$, or the $j$-th source and mixture, respectively.

\section{The Convex Cone Method}
\subsection{NN Approach}
\begin{figure}
\includegraphics[height=5cm,width=8cm]{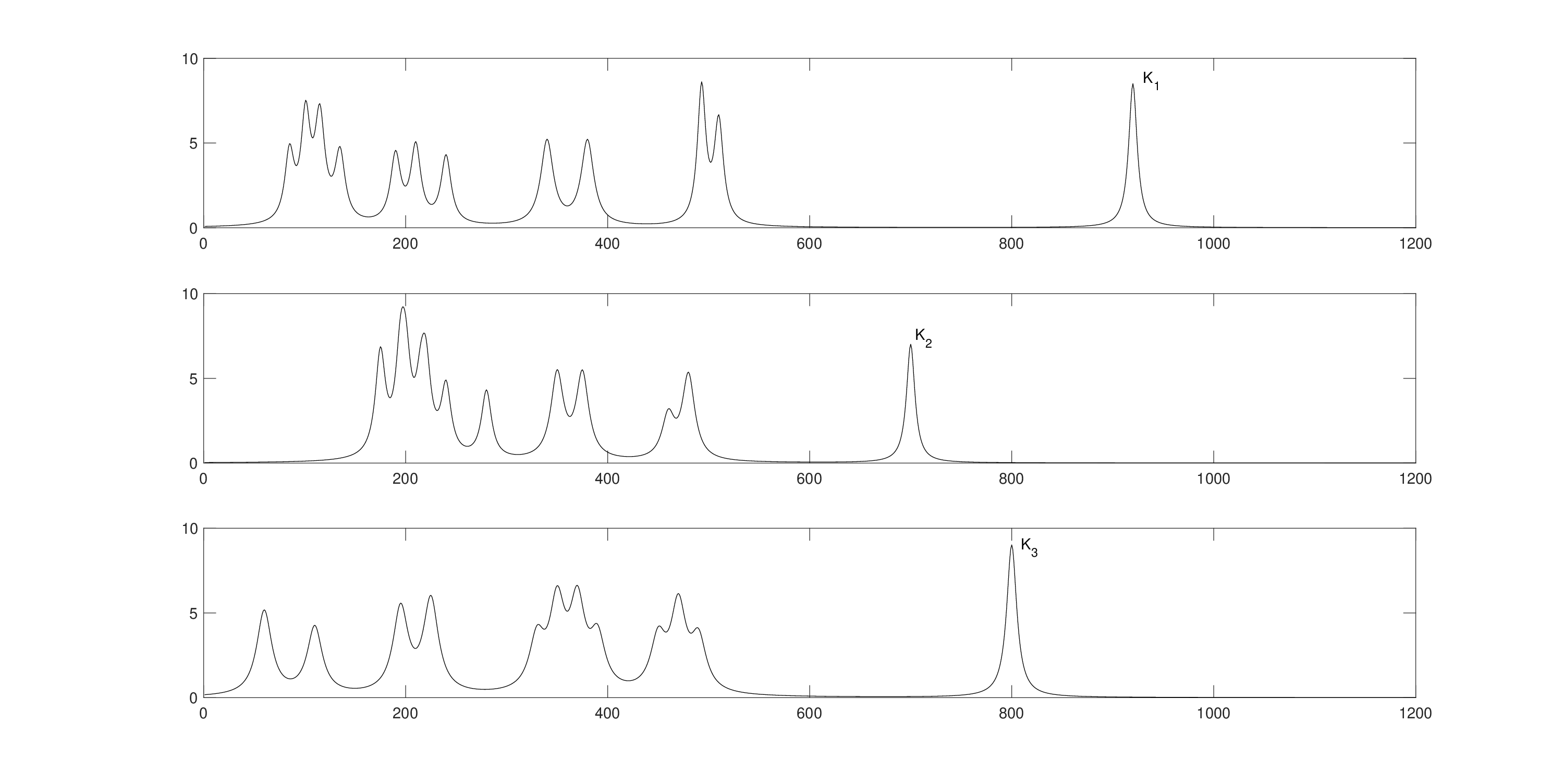}
\includegraphics[height=5cm,width=8cm]{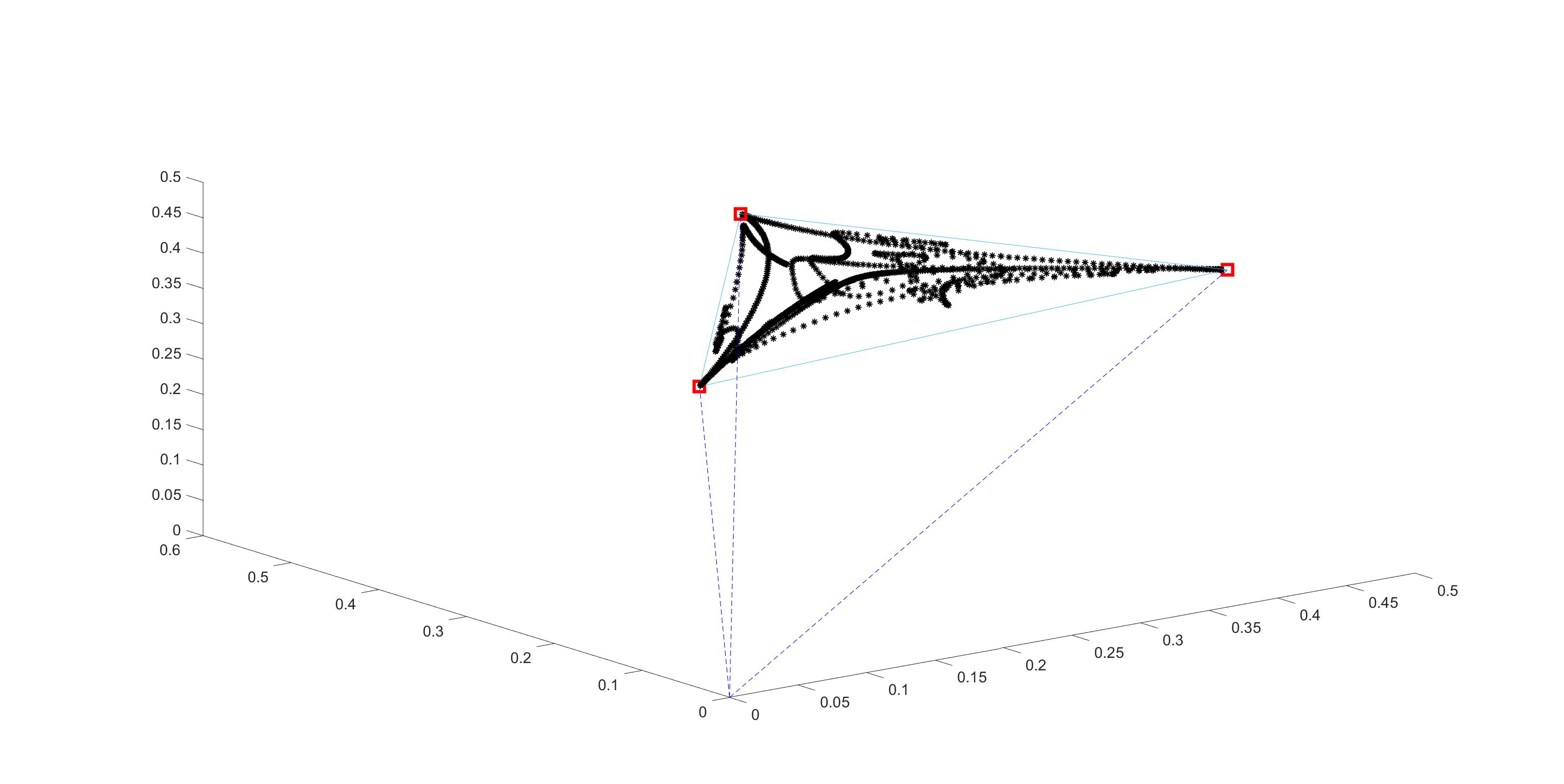}
\caption{ Left : simulated NMR spectra of three SAP sources.  Each signal has a stand alone peak indicated by $K_1,K_2$, and $K_3$.  Right: the scattered plot of $X$ (columns of $X$) scaling to be on plane $x+y+z = 1 $.}
 \label{NNA}
\end{figure}

In \cite{NN05}, Naanaa and Nuzillard (NN) presented an efficient sparse BSS method and its mathematical analysis for nonnegative and partially orthogonal signals such as NMR spectra.  Consider the determined regime where there are same number of mixtures is no less than that of sources ($m\geq n$), and the mixing matrix $A$ is full rank.   In simple terms, NN's key sparseness assumption (referred to as NNA below) on source signals is that each source has a stand-alone peak at some
location of the acquisition variable where the other sources are identically zero.  More precisely, the source matrix $S\geq 0$ is assumed to satisfy the following condition
\begin{ass}[NNA]: For each $i\in\{1,2,\dots,n \}$ there exists an $j_{i}\in
\{1,2,\dots,p\}$ such that $s_{i,j_i}>0$ and $s_{k,j_i}=0\; (k=1,\dots,i-1,i+1,\dots,n)\;.$
\end{ass}
Eq. (\ref{LMM}) can be rewritten in terms of columns as
\begin{equation}
\label{LinComb}
 X^{j} = \sum^{n}_{k=1}s_{k,j}A^{k},\;\; j = 1,\dots,p\;,
\end{equation} where $ X^j$ denote the $j$th column of $X$, and $A^k$ the $k$th column of $A$.
Assumption NNA implies that $\displaystyle X^{j_i} = s_{i,j_i}A^i\;\; i = 1,\dots,n\;\; $ or $A^{i} = \frac{1}{s_{i,j_i}}X^{j_i}$. Hence
Eq. (\ref{LinComb}) is rewritten as
\begin{equation}
\label{NNlinComb} X^{j} = \sum^{n}_{i = 1} \frac{s_{i,j}}{s_{i,j_i}}X^{j_i}\;,
\end{equation}
which says that every column of $X$ is a nonnegative linear combination of the columns of $\hat{A}$.  Here $\hat{A} =
[X^{j_1},\dots,X^{j_n}]$ is the submatrix of $X$ consisting of $n$ columns each of which is collinear to a particular column of
$A$.  It should be noted that $j_i\;(i = 1,\dots,n) $ are not known and have to be computed.  Once all the $j_i s$ are found,
an estimation of the mixing matrix is obtained.  The identification of $\hat{A}'s$ columns is equivalent to identifying a convex cone of a finite collection of vectors.  The cone encloses the data columns in matrix $X$, and is the smallest of such cones.  Such a minimal enclosing convex cone can be found by linear programming methods.  Mathematically, the following constrained equations are formulated for the identification of $\hat{A}$,
\begin{equation}
\label{LPNF} \sum^{p}_{j = 1, j\neq k}X^{j} \lambda_j = X^{k}\;, \lambda_j\geq 0\;\; k = 1,\dots,p\;.
\end{equation}
Then any column $X^{k}$ will be a column of $\hat{A}$ if and only if the constrained equation (\ref{LPNF}) is inconsistent.
However, if noises are present, the following optimization problems are suggested to estimate the mixing matrix
{\allowdisplaybreaks
\begin{eqnarray*}
\label{LPNP} & &\mathrm{minimize}\;\; \mathrm{score} = \|\sum^{p}_{j = 1, j\neq k}X^{j} \lambda_j - X^{k}
 \|_2\;, k = 1,\dots,p\;\\
 & & \mathrm{subject \; to} \;\; \lambda_j \geq 0\;.
\end{eqnarray*}}
A score is associated with each column.  A column with a low score is unlikely to be a column of $\hat{A}$ because this
column is roughly a nonnegative linear combination of the other columns of $X$.  On the other hand, a high score means that the corresponding column is far from being a nonnegative linear combination of other columns.  Practically, the $n$ columns from $X$ with highest scores are selected to form $\hat{A}$, the mixing matrix.  The Moore-Penrose inverse $\hat{A}^{+}$ of $\hat{A}$ is then computed and an estimate to $S$ is obtained: $\hat{S} = \hat{A}^{+} X$.   NN method proves to be both accurate and efficient if NNA condition holds.  However, if the condition is not satisfied, errors and artifacts may be introduced because the true mixing matrix is no longer the smallest enclosing convex cone of columns of the data matrix. As it applies to signals with peak structures such NMR spectra, NNA can be restated as the stand alone peak (SAP) condition: each source signal possesses a stand alone peak over certain acquisition interval, where other sources are identically zero. SAP condition is illustrated by NMR spectra of three sources in the left plot of Fig. \ref{NNA}, it can be seen that each source signal has a stand alone peak denoted by $P_1,P_2,$ and $P_3$, respectively.

Recently, the authors have developed postprocessing techniques on how to improve NN results with abundance of mixture data, and how to improve mixing matrix estimation with major peak based corrections \cite{sun_xin_pNN}.
The work actually considered a relaxed NNA (rNNA) condition
\begin{ass}[rNNA]: For each
$i\in\{1,2,\dots,n \}$ there exists an $j_{i}\in \{1,2,\dots,p\}$ such that $s_{i,j_i}>0$ and $s_{k,j_i}\approx\epsilon_k\;
(k=1,\dots,i-1,i+1,\dots,n)\;,$ where $ s_{i,j_i}\gg \epsilon_k$.
 \end{ass}
Simply said, each source signal has a dominant peak at acquisition position where the other sources are allowed to be nonzero. NNA condition recovers if all $\epsilon_k = 0$.  The rNNA is more realistic and robust than the ideal NNA for real-world NMR data \cite{NN05}.
\subsection{Dominant Intervals and Degenerate Matrix}
\begin{figure}
\includegraphics[height=5cm,width=8cm]{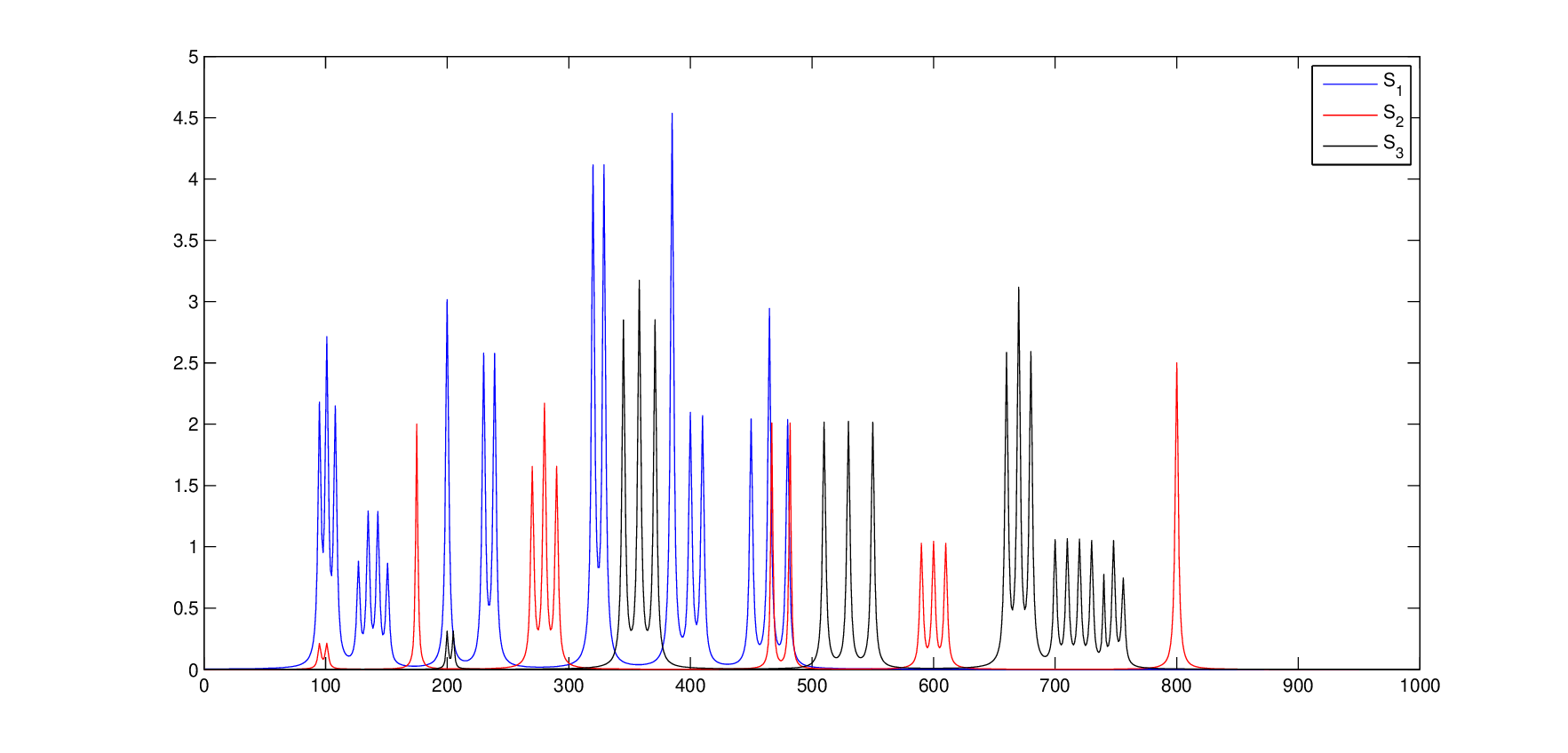}
\includegraphics[height=5cm,width=8cm]{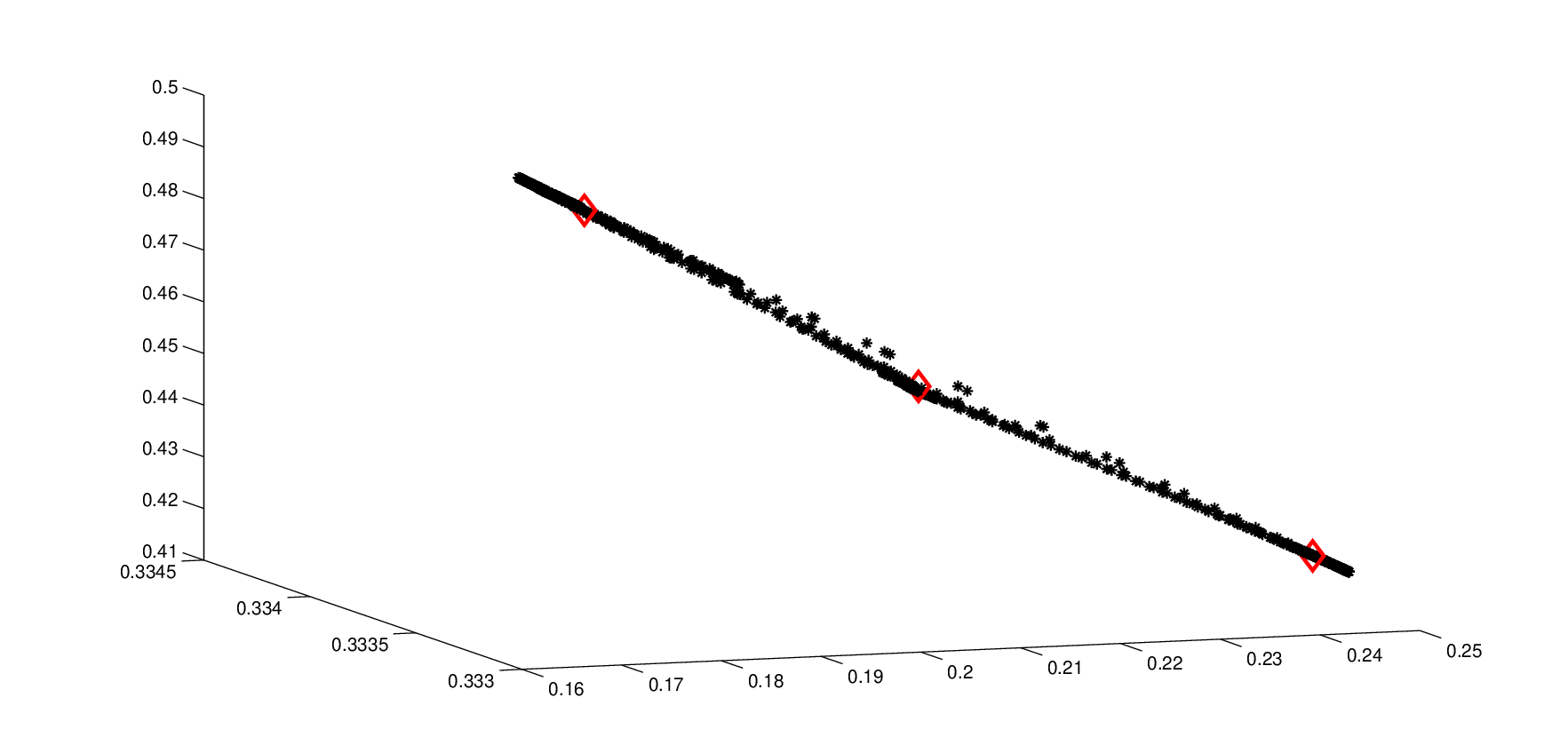}
\caption{three source signals with dominant intervals (left panel); the geometry of the mixture matrix (right panel).  The centers (red diamond) of three clusters are detected by k-means. }
\label{SourceCondition}
\end{figure}

Motivated by the DOSY NMR spectra, we propose here a different relaxed NN condition on the source signals.  Note that the rows $S_1,S_2,\dots,S_n$ of $S$ are the source signals, and they are required to satisfy the following condition:  For  $i =1,2,3,\dots,n$, source signal $S_i$ is required to have dominant interval(s) over $S_{n},\dots,S_{i+1}, S_{i-1},\dots,S_2,S_1$, while $S_i$ is allowed to overlap with other signals at the rest of the acquisition region.  More formally, this condition implies that source matrix $S$ satisfies the following condition
\begin{ass}
For each $k \in \{1,2,3,\dots,n \}$,  there is a set $\mathcal{I}_k \subset \{1,2,\dots,p \}$ such that for each $l\in \mathcal{I}_k$ $s_{k l} \gg s_{j l}, j = 1,2,\dots,k-1,k+1,\dots, n$.
\end{ass}
We shall call this dominant interval condition, or DI condition.
Fig. \ref{SourceCondition} is an idealized example of three DI source signals. The motivation lies in the similar diffusion rates of the chemicals with similar structure.  This poses a mathematical challenge to invert a nearly singular matrix, since a small error in the recovered mixing matrix might lead to a considerable deviation in the source recovery.  Among the singularly mixed signals (or degenerate data), in this paper we shall consider the following two types: 1) columns of the mixing matrix are parallel; 2) one column of the mixing matrix is a nonnegative linear combination of others.  Case 1 is  motivated by NMR of the chemicals with similar diffusion rates.  We shall call this condition parallel column condition, or PCC.  Case 2 can also be encountered in NMR spectroscopy of chemicals, and we shall call it one column degenerate condition, or OCDC.  Throughout the paper, we shall refer to a mixing matrix of PCC or OCDC as degenerate mixing matrix. Please note that both PCC and OCDC should be considered to hold approximately in real-world data.
\section{Our Methods}
\subsection{Data Clustering}
\begin{figure}
\includegraphics[height=5cm,width=8cm]{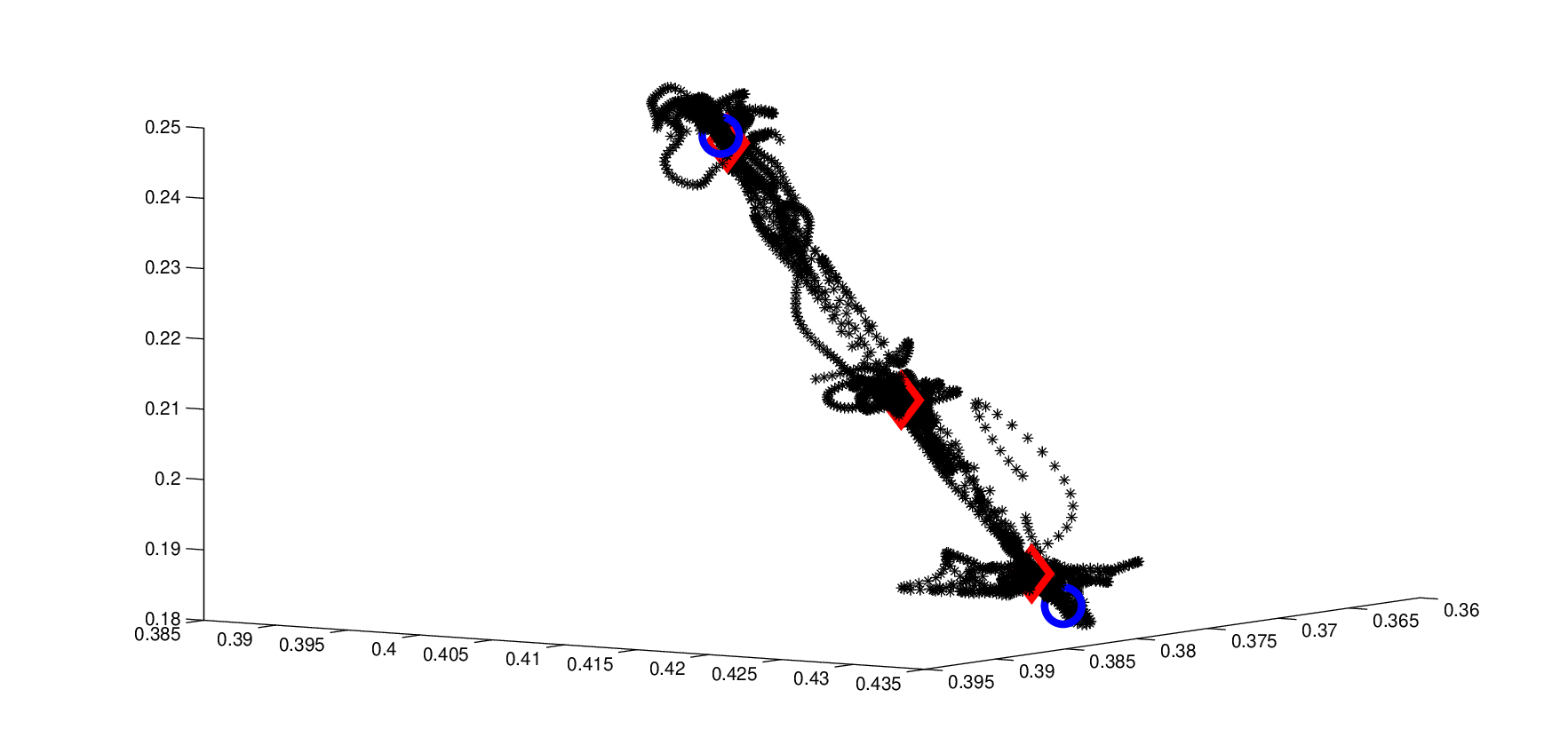}
\includegraphics[height=5cm,width=8cm]{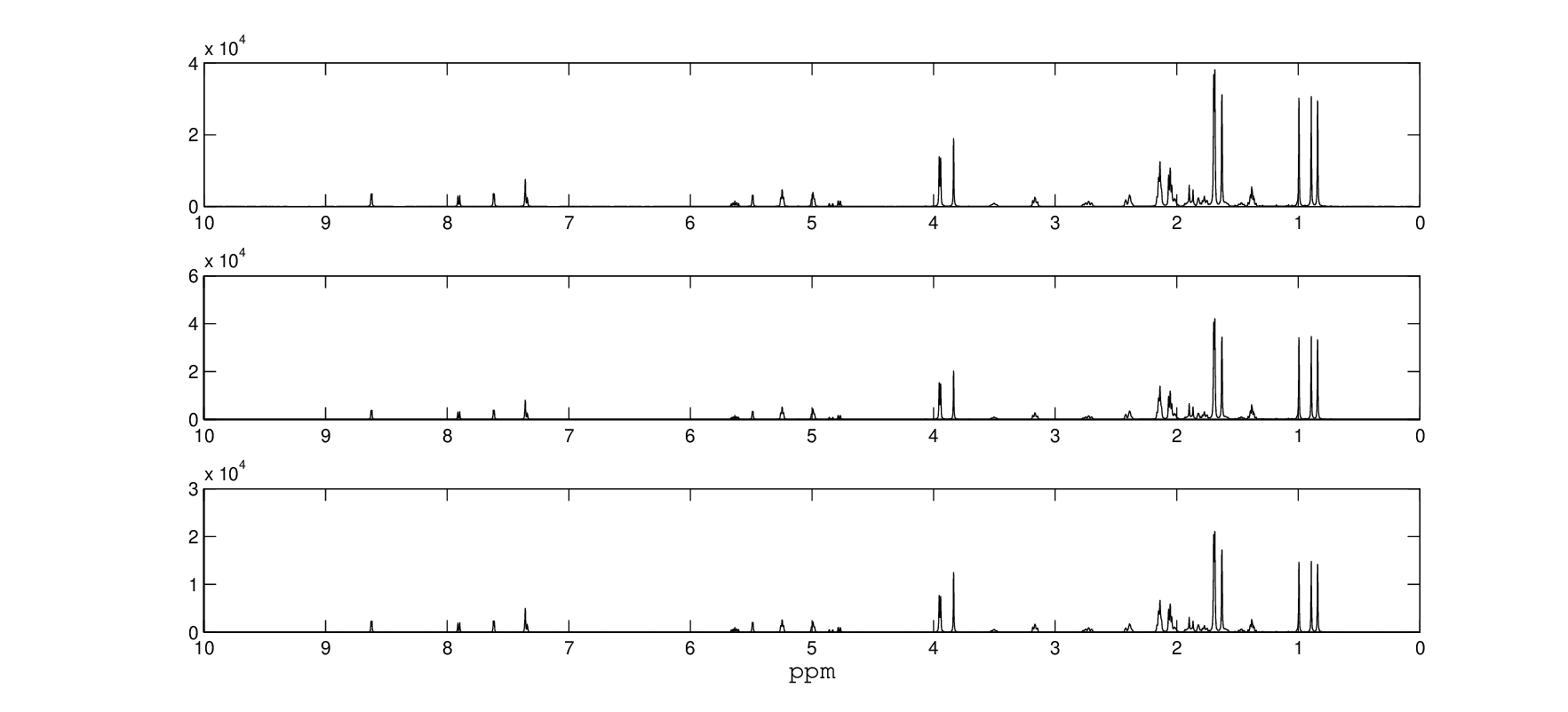}
\caption{Real data example: three columns of $A$ are identified as the three center points (in red diamond) attracting most points in scatter plots of the columns of X (left), and the three rows of X (right).  NN method identifies two columns of $A$ as the points in the blue circle.}
\label{real_data}
\end{figure}

\begin{figure}
\includegraphics[height=5cm,width=8cm]{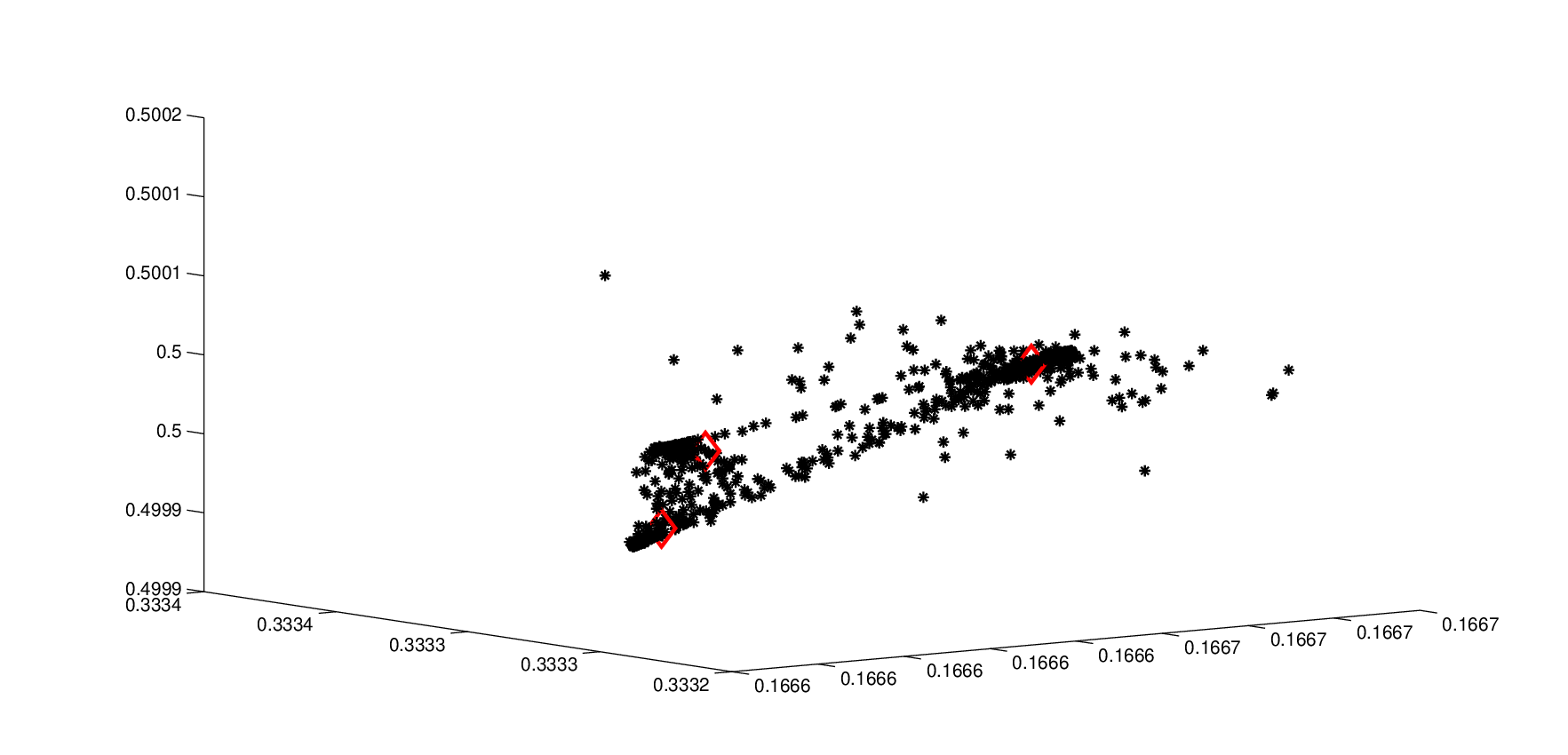}
\includegraphics[height=5cm,width=8cm]{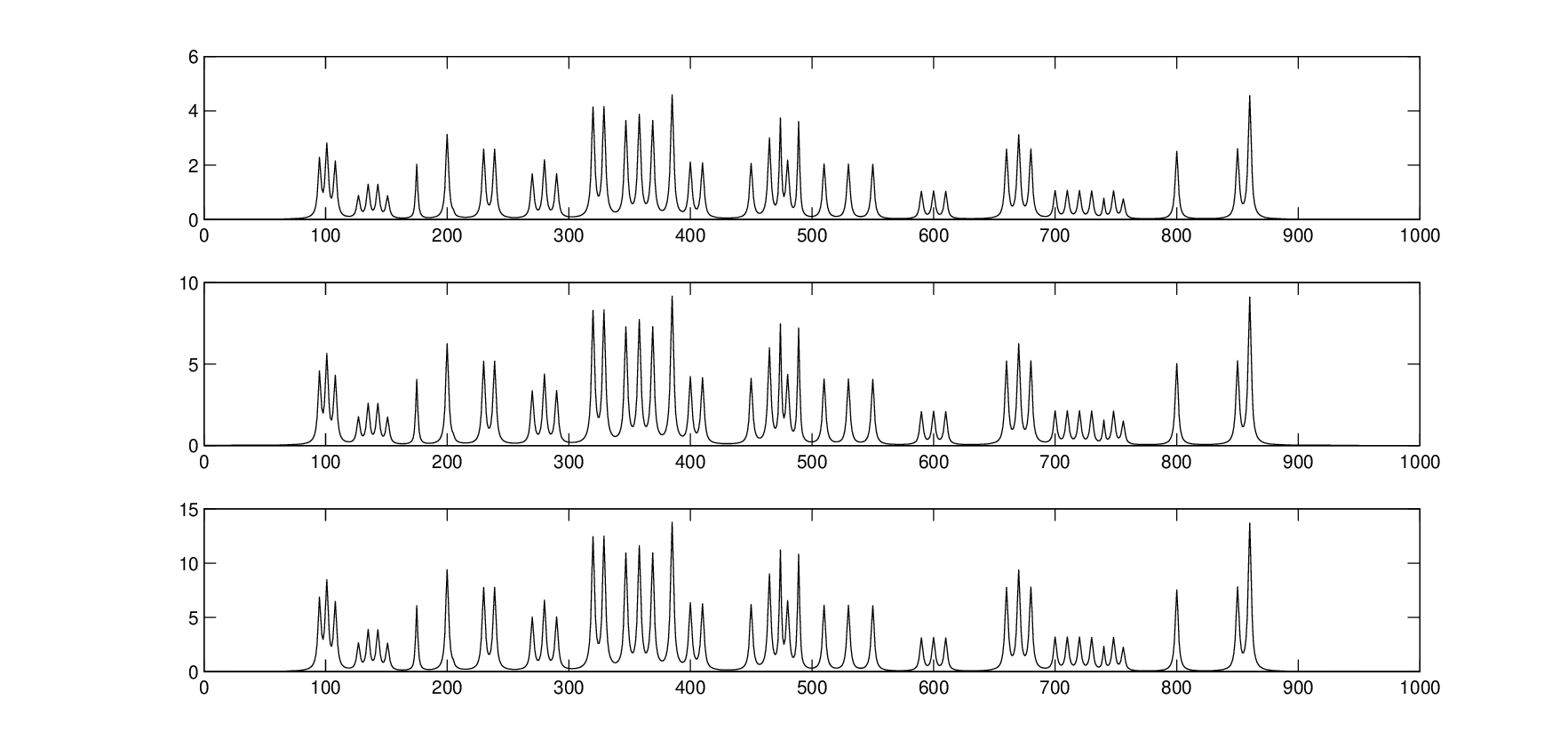}
\caption{Example of PCC case: three columns of $A$ are identified as the three center points (in red diamond) attracting most points in scatter plots of the columns of X (left), and the three rows of X (right).}
\label{real_data_PCC}
\end{figure}

Now suppose we have a set of nearly degenerate signals from DI sources.  We require that compared to the size of dominant interval(s) in the acquisition region, the source signals overlapping region is much smaller.  In fact, this is a reasonable assumption for the NMR data which motivates us.  More importantly, this requirement facilitates the successful implementation of the clustering method.  Next, we shall estimate the columns of mixing matrix $A$ by clustering.  The dominant interval(s) from each of the source signals implies that there is a region where the source $S_\mathbf{i}$ dominates others.  More precisely, there are columns of $X$ such that
\begin{equation}
\label{multi_X}
X^k = s_{{i},k}A^i + \sum^{n-1}_{j = 1}o_{j,k}A^j\;,
\end{equation}
where $s_{{i},k}$ dominate $o_{i,k} (i = 1,\dots,n-1)$, i.e., $s_{{i},k} \gg o_{i,k} $.  The identification of $A^i$ $(i = 1,\dots, n)$ is equivalent to finding a cluster formed by these $X^k$'s in $\mathbb{R}^m$.  As illustrated in the geometry plot of $X$ in Fig. \ref{SourceCondition}, three clusters are formed.  Many clustering techniques are available for locating these clusters, for example, k-means is one of the simplest unsupervised learning algorithms that solve the well known clustering problem.  We shall use k-means for the data in this paper, returns the centers of the clusters.  Those centers are the estimates of columns of $A$.  Consider an example of three DI source signals with OCDC mixing matrix condition, the three centers shown in Fig. \ref{SourceCondition}.  For real-world data, we show an example of NMR spectra of quinine, geraniol, and camphor mixture in Fig. \ref{real_data}.  The clusters in the middle implies that OCDC condition hold well for this data.  Apparently, NN method (and other convex cone methods) would fail to separate the source signals due to the degeneracy of the mixing matrix.  It might be able to identify two columns of $A$ as the two edges, it by no means can locate $A$'s degenerate column.  For the PCC degenerate case, clustering is also able to deliver a good estimation, even when the data is contaminated by noise.  We show the results in Fig. \ref{real_data_PCC} where the three clusters are very close due to the PCC degeneracy.  NN solution would deviate considerably from the true solution.  For the data we tested, clustering techniques like k-means works well when the condition number of the mixing matrix is up to $10^8$.  Though the solutions of mixing matrix by clustering methods are rather good estimation to the true solution, small deviations from the true ones will introduce large errors in the source recovery ($S = \mathrm{inverse}(A)\,X$).  Next we propose two approaches to overcome this difficulty.  Both approaches need to solve optimization problems. The first one intends to improve the source recovery by seeking a better mixing matrix, while the second approach reduces the spurious peaks by imposing sparsity constraint on the sources.

\subsection{Better Inverse of the Mixing Matrix}
Suppose the estimation of the mixing matrix by clustering is $\hat{A}$.  As discussed above, errors in $\hat{S}$ could be introduced even by a small deviation in $\hat{A}$ from the ground truth. Negative spurious peaks are produced in most cases, see the Fig. \ref{eg1_rec} where the negative peaks on the left plot actually can be viewed as bleed through from another source.  Clearly, a better estimation of mixing matrix is required to reduce these spurious peaks.  Instead of looking for a better mixing matrix, we propose to solve the following optimization problem for a better inverse of the matrix in the determined case ($m=n$, same number of mixtures of sources),
\begin{equation}
\label{inverse}
 \min_{B} \frac{1}{2}\|I-B\hat{A}\|^2_2\quad \mathrm{subject}\;\mathrm{to}\; B\,X\geq 0\;,
\end{equation}
where $I\in \mathbb{R}^{n\times n}$ is the identity matrix.  The constraint $B\,X \geq 0$ is used to reduce the negative values introduced in the source recovery.  (\ref{inverse}) is a linearly constrained quadratic program and it can be solved by a variety of methods including interior point, gradient projection, active sets, etc.  In this paper, interior point algorithm is used.  Once the minimizer $B^*$ is obtained, we solve for the sources by $S = B^*\,X$.  For over-determined mixture ($m>n$, $A$ being a tall matrix), the source recovery is usually obtained by $\hat{S} = (\hat{A}^T\hat{A})^{-1}\hat{A}^T X$, to achieve a better inverse of $\hat{A}$ in the sense of least squares, the following optimization problem is proposed,
\begin{equation}
\label{inverseLS}
 \min_{B} \frac{1}{2}\|I-B\hat{A}^T\hat{A}\|^2_2\quad \mathrm{subject}\;\mathrm{to}\; B\,\hat{A}^T X\geq 0\;,
\end{equation}
which can be solved via an interior point algorithm as well.
\subsection{Sparser Source Signals}
The method proposed above works well for mixing matrix whose condition number is up to $10^8$.  If the mixing matrix is much more ill-conditioned, the problem (\ref{LMM}) becomes under-determined.  It appears that solving the equation exactly for $S$ is hopeless even an accurate $A$ is provided.  However, a meaningful solution is possible if the actual source signals are structurally compressible, meaning that they essentially depend on a low number of degrees of freedom.  Although the source signals (rows of $S$) are not sparse, the columns of $S$ possess sparsity due to the dominant intervals condition.  Hence, we seek the sparsest solution for each column $S^i$ of $S$ as
\begin{equation}
\label{Lzero}
 \min \|S^i \|_0\quad \mathrm{subject}\;\mathrm{to}\; \hat{A}S^i = X^i,\;\; S^i \geq 0.
\end{equation}
Here $\|\cdot\|_0$ ( 0-norm ) represents the number of nonzeros.
Because of the non-convexity of the 0-norm, we minimize the $\ell_1$-norm:
\begin{equation}
 \label{Lone}
 \min \|S^i \|_1\quad \mathrm{subject}\;\mathrm{to}\; AS^i = X^i,\;\; S^i \geq 0,
\end{equation}
which is a linear program \cite{Don} because $S^i$ is non-negative.
The fact that data may in general contain noise suggest solving the following unconstrained optimization problem,
\begin{equation}
 \label{LoneU}
 \min_{S^i \geq 0} \mu\|S^i \|_1 + \frac{1}{2}\|X^i - AS^i \|^2_2\;,
\end{equation}
for which Bregman iterative method \cite{G_O_09, YO} with a proper projection onto non-negative convex subset
can be used to obtain a solution.  Under certain conditions of matrix $A$, it is known \cite{CanT,Zhang}
that solution of $\ell_1$-minimization (\ref{LoneU}) gives the
exact recovery of sufficiently sparse signal, or solution to
(\ref{Lzero}), \cite{CanT,Zhang}.  Though our numerical results support the
equivalence of $\ell_1$ and $\ell_0$ minimizations, the mixing matrix $A$ does not satisfy the existing
sufficient conditions \cite{CanT,Zhang}.

\subsection{Structure Assisted Nonnegative Matrix Factorization}
In the case where the source signals are known to possess dominant intervals and clusters in the data points, convex BSS methods proposed above are suitable choices.  While if no knowledge of the source signals is available other than the nonnegativity, we shall opt to use nonnegative matrix factorization approach.
Nonnegative matrix factorization (NMF) \cite{Lee,NMF0} is the prevalent method for solving nonnegative BSS problems with limited information on the nonnegative source signals and mixing matrices.  As a parts-based data representation, NMF attemptes to find two nonnegative matrices whose product approximates the data matrix.  Given the mixture matrix $X \in \mathbb{R}^{m\times p}$, NMF seeks the approximate factorization  $X \approx A S$, $A\in \mathbb{R}^{m\times n}, S \in \mathbb{R}^{n\times p} $ by minimizing the following objective function,
\begin{equation*}
J(A,S) = \|\ X- A\,S\|_{F}^2\;,
\end{equation*}Where $\|\cdot\|_F$ is the Frobenius matrix norm.  $J(A,S)$ is also called Eucleadian distance function which is non-convex in $A,S$, so it is unrealistic to expect an algorithm to converge to the global minimum.  Lee and Seung \cite{Lee} proposed a multiplicative iteration rule:
\begin{eqnarray}
  a_{ij} &\leftarrow& a_{ij}\frac{\left[ XS^T \right]_{ij}}{\left [ ASS^T \right]_{ij}},\\
  s_{jk} &\leftarrow& s_{jk}\frac{\left[ A^TX\right]_{jk}}{\left [ A^TAS \right]_{jk}};,
\end{eqnarray} and they proved that these iterates converge to a local minimum of $J(A,S)$.

Given the fact that the columns of the mixing are similar (PCC columns). i.e., the sources signals are accounts for approximately same amounts.  We should incorporate this constraint into the $J(A,S)$ to steer the solution towards the desired one. The linear mixture model $X = AS$ can be rewritten as follows:
\begin{equation}
\label{bssColumn}
[X^1,X^2,\cdots,X^p] = [A^1,A^2,\cdots,A^n]\left(
\begin{array}{ccccc}
     S_{11} & S_{12}   & \cdots  & S_{1p}\\
     S_{21} & S_{22}   & \cdots  & S_{2p} \\
      \vdots & \vdots  & \vdots  & \vdots \\
     S_{n-1,1} & S_{n-1,2}& \cdots & S_{n-1,p}\\
     S_{n1} & S_{n2}& \cdots & S_{n,p}
       \end{array}
 \right)\;,
\end{equation}
We assume that matrix $A$ is nearly singular as we assume that its columns $A^1,A^2,\dots, A^n$ are similar. This situation holds approximately well in mixtures where the source signals have nearly the same amounts of concentration; or in image processing such as computer vision, the adjacent frames should vary very little, or look similar.  With limited knowledge on the source signals, we shall opt for nonnegative matrix factorization approach \cite{Lee}, we need to build the mixing matrix structure into the model.  The similarity between the columns of $A$ implies that the differences of them (given that the columns are all normalized to unit vectors) are approximately zeros, or
\begin{equation*}
A^2-A^1 \approx \mathbf{0}, A^3-A^2 \approx \mathbf{0}, \cdots, A^{i}-A^{i-1} \approx \mathbf{0}, \cdots A^{n}-A^{n-1}\approx \mathbf{0}, A^1-A^{n} \approx \mathbf{0}.
\end{equation*}
To account for these similar (parallel) columns, we shall put these conditions into matrix form
\[
\begin{bmatrix}
A^1 & A^2 & A^3 & \cdots & A^{n-1} & A^{n}\\
\end{bmatrix}
\begin{bmatrix}
-1     & 0    & 0       & \dots      &0        & 1   \\
1      & -1   & 0       & \dots      & 0       & 0   \\
0      & 1    & -1      & \dots      & 0       & 0  \\
\vdots &\vdots & \vdots &  \vdots    &\vdots   & \vdots  \\
0       & 0    & 0      & \dots       & -1      & 0\\
0       & 0    & 0       & \dots       & 1  & -1 \\
\end{bmatrix}
\approx
\mathbf{O}
\]
Where denote $W \in \mathbb{R}^{n\times n}$ as the sparse matrix containing $\pm 1 's$.  $\mathbf{O}\in \mathbb{R}^{m\times n} $ is the zero matrix. Now let us consider the following cost function
\begin{equation}
\label{sNMF}
J_{\alpha}(A,S) = \|\ X- A\,S\|_{F}^2 + \alpha \|A W\|_F^2, \mathrm{s.t. } A \succeq 0, S \succeq 0.
\end{equation}
  Following \cite{Choi} by the gradient descent, we are able to get the update rules for entries of  $A$ and $S$
\begin{eqnarray*}
  a^{l+1}_{ij} - a^{l}_{ij} & = &  -\delta_{ij} \frac{\partial J}{\partial a_{ij} }\\
   & = & \delta_{ij} \left ( [XS^{T}]_{ij} - [ASS^T]_{ij}-\alpha [AWW^T]_{ij} \right ) \\
  s^{l+1}_{jk} - s^{l}_{jk} & = & -\eta_{jk} \frac{\partial J}{\partial s_{jk} }\\
   & = & \eta_{jk} \left ( [AX^T]_{jk} - [A^TAS]_{ij} \right ) \\
\end{eqnarray*}
To ensure the nonnegativity, we shall choose the specific learning rates proposed in \cite{Lee}
$$
\delta_{ij} = \frac{a_{ij}}{\left [ ASS^T \right]_{ij}}, \eta_{jk} = \frac{s_{jk}}{\left [ A^TAS \right]_{jk}}
$$
 and the resulting multiplicative update formulas are
\begin{eqnarray}
  a_{ij} &\leftarrow&  \frac{a_{ij}}{\sum_i a_{ij}}, \\
  a_{ij} &\leftarrow& a_{ij}\frac{\left[  (XS^T-\alpha AWW^T)_{ij},0\right]_{+}}{\left [ ASS^T \right]_{ij}},\\
  s_{jk} &\leftarrow& s_{jk}\frac{\left[ A^TX\right]_{jk}}{\left [ A^TAS \right]_{jk}};,
\end{eqnarray} where the operator $[x]_{+} = \max{0, x}$ is to ensure nonnegativity in the updates.  The proposed NMF cost function addresses the degeneracy of the mixing matrix by imposing a penalty on the differences of its columns, and shall be termed by column difference NMF, or \textbf{CD-NMF}.

Below we propose an alternative cost function for NMF.  Noticing that to measure the similarity between vectors, we can use the inner product as a simple yet accurate measure for similarity.  For two parallel unit vectors $a, b$, their inner product $a\cdot b = 1$. Based on this idea, we propose to adopt the following criterion to impose the similarity between the columns of the mixing matrix $A$.  Then we shall have $A^{T}A = J_n$, here $J_n$ is a square matrix of size $n$ of ones, that is,
\[
A^{T}A =
\begin{bmatrix}
{A^1\cdot A^1}     & A^1\cdot A^2             & A^1\cdot A^3       & \dots      & A^1\cdot A^n \\
        A^2\cdot A^1      & {A^2\cdot A^2}    & A^2\cdot A^3       & \dots      & A^2\cdot A^n          \\
           \vdots          &\vdots                    & \vdots            &  \vdots    &  \vdots     \\
        A^n\cdot A^1       & A^n\cdot A^2             &A^n\cdot A^3       & \dots      & {A^n\cdot A^n}\\
\end{bmatrix}
= \begin{bmatrix}
1     & 1             & 1       & \dots      & 1 \\
        1     & 1    & 1       & \dots      & 1         \\
           \vdots          &\vdots                    & \vdots            &  \vdots    &  \vdots     \\
       1       & 1             &1      & \dots      & 1\\

\end{bmatrix}= J_n
\]
where we can see that the normalization of columns of $A$, or $A^1, A^2,\cdots, A^n$ are all enforced as the diagonals show. We consider the following cost function
\begin{equation}
\label{pNMF}
J_{\beta}(A,S) = \|\ X- A\,S\|_{F}^2 + \beta \|A^TA - J_n\|_F^2, \mathrm{s.t. } A \succeq 0, S \succeq 0.
\end{equation}
By a likewise derivation presented in the previous section, we have
\begin{eqnarray*}
  a^{l+1}_{ij} - a^{l}_{ij} & = &  -\delta_{ij} \frac{\partial J}{\partial a_{ij} }\\
   & = & \delta_{ij} \left ( [XS^{T}]_{ij} - [ASS^T]_{ij}-2\beta [AA^TA-AJ_n^T]_{ij} \right ) \\
  s^{l+1}_{jk} - s^{l}_{jk} & = & -\eta_{jk} \frac{\partial J}{\partial s_{jk} }\\
   & = & \eta_{jk} \left ( [AX^T]_{jk} - [A^TAS]_{jk} \right ). \\
\end{eqnarray*}
 We again take the previous learning rates
 $$
\delta_{ij} = \frac{a_{ij}}{\left [ ASS^T \right]_{ij}}, \eta_{jk} = \frac{s_{jk}}{\left [ A^TAS \right]_{jk}}
 $$
which lead us to the following update rules:
\begin{eqnarray}
  a_{ij} & \leftarrow & a_{ij}\frac{\left[ [XS^T]_{ij}-2\beta[AA^TA-AJ_n^T]_{ij},0 \right ]_{+}}{\left [ ASS^T \right]_{ij}},\\
  s_{jk} & \leftarrow & s_{jk}\frac{\left[ A^TX\right]_{jk}}{\left [ A^TAS \right]_{ij}};\\
\end{eqnarray}
This proposed NMF cost function addresses the degeneracy of the mixing matrix by adding a penalty on the inner products  of all its columns, and shall be termed by column product NMF, or \textbf{CP-NMF}.

\section{Numerical experiments}

\subsection{Singular Mixing Matrix and Dominant Peak Sources}
\begin{figure}
\includegraphics[height=5cm,width=8cm]{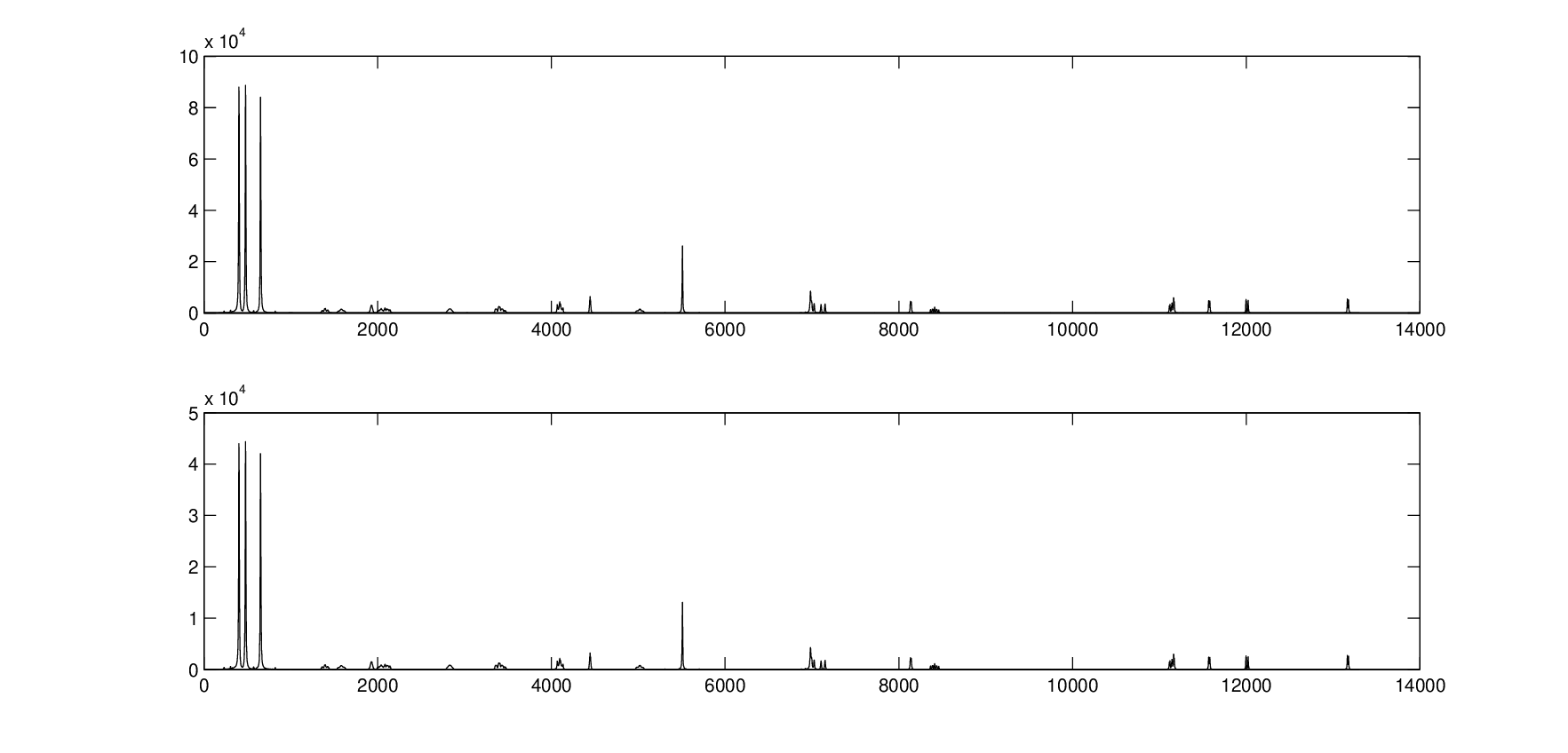}
\includegraphics[height=5cm,width=8cm]{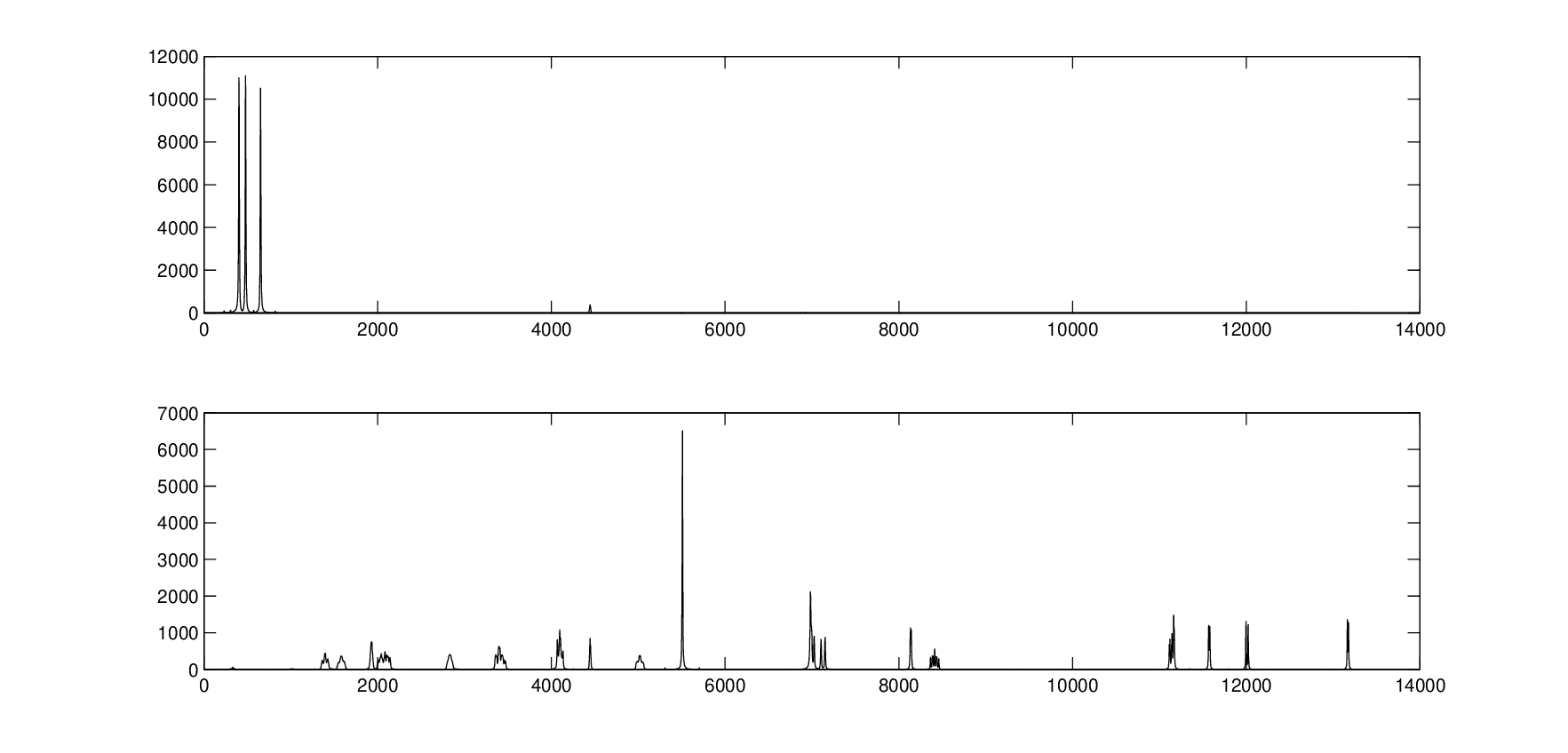}
\caption{recovered sources by clustering (left column) and the ground truth (right column).}
\label{eg1_mix}
\end{figure}

\begin{figure}
\includegraphics[height=5cm,width=5cm]{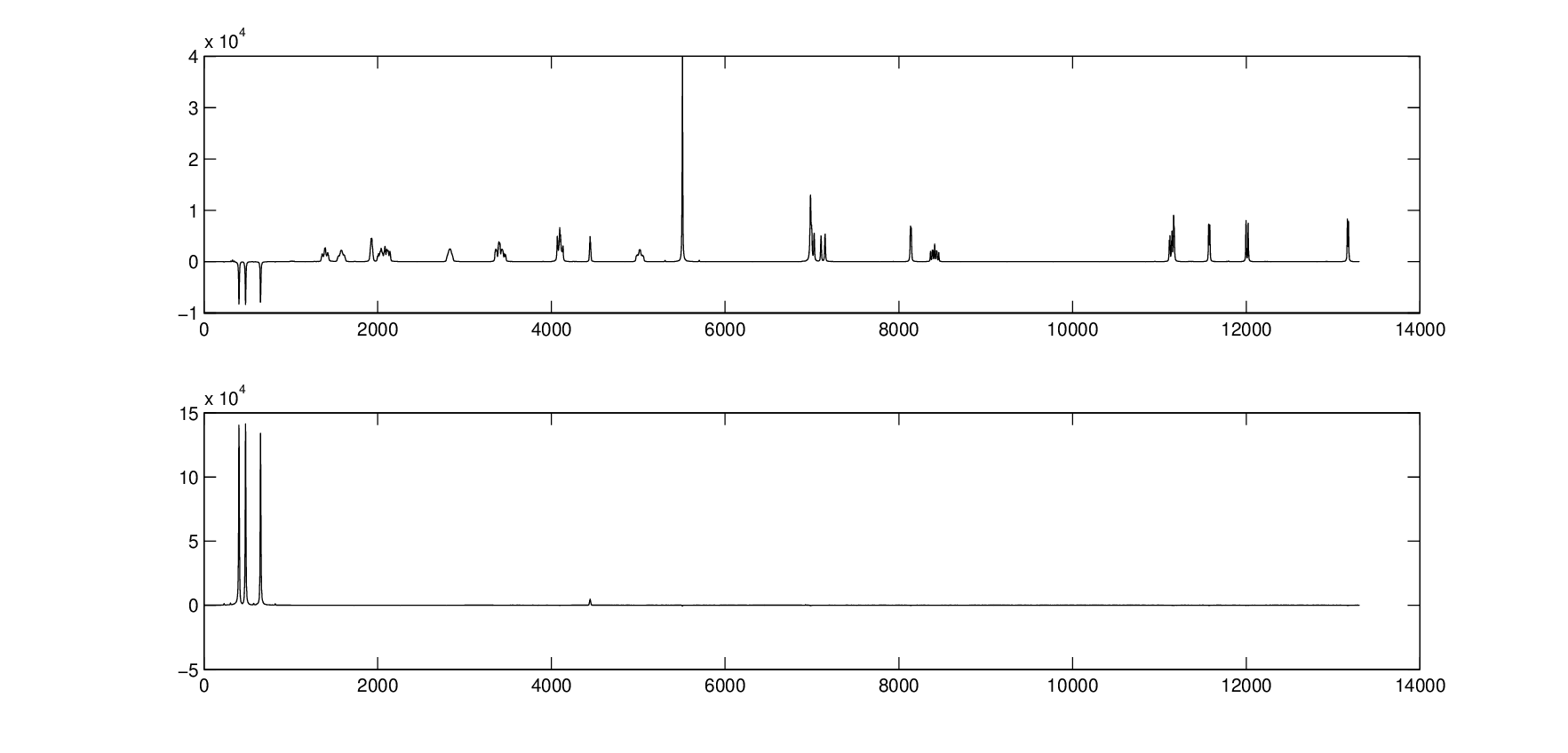}
\includegraphics[height=5cm,width=5cm]{figure/source_ref.eps}
\includegraphics[height=5cm,width=5cm]{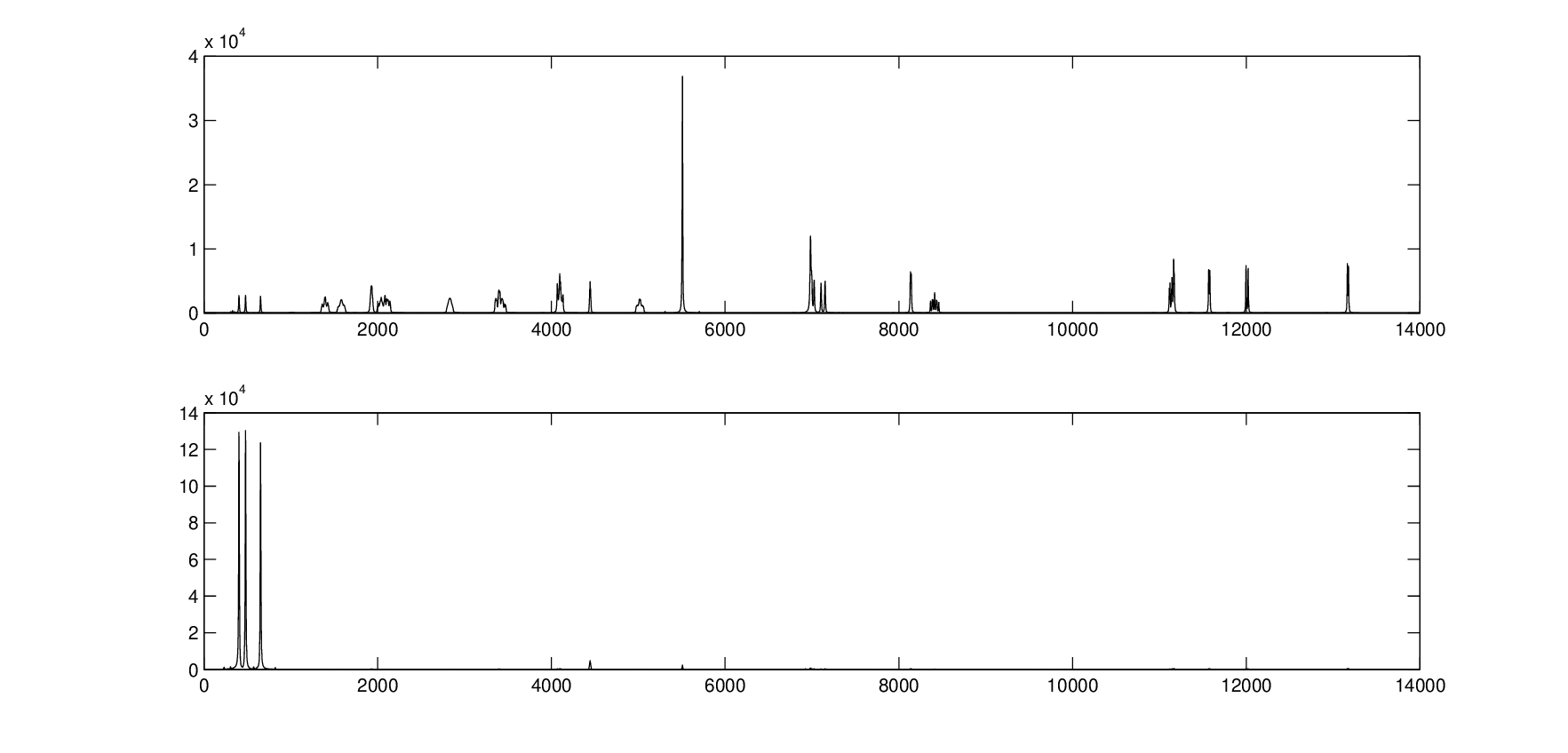}
\caption{recovered sources by clustering (left column), the ground truth (middle column), and the improved results by a better estimate of mixing matrix (right column).}
\label{eg1_rec}
\end{figure}

In this section, we report the numerical examples solved by the method.   We compute three examples.
The data of the first two examples are synthetic, while the third example uses real NMR data.
In the first example, two sources are to be separated from two mixtures.  The mixtures are constructed from two real NMR source signal by simulating the linear model (\ref{LMM}).   The two columns of mixing matrix are nearly parallel, and its condition number is about $1.25\times 10^8$.  The true mixing matrix $A_{\mathbf{TR}}$, its estimation $A_\mathbf{C}$ via clustering, and the improved estimate $A_\mathbf{P}$ by solving (\ref{inverse}) are
(for ease of comparison, the first rows of $\hat{A},A_p$ are scaled to be same as that of $A$)
\begin{equation*}
  A_{\mathbf{TR}} =\left(
   \begin{array}{ccccc}
    0.894427190999916  & 0.894427182055644\\
    0.447213595499958 &  0.447213613388501
   \end{array}
 \right)\;,
\end{equation*}
\begin{equation*}
  A_\mathbf{C} =\left(
   \begin{array}{ccccc}
    0.894427190999916 & 0.894427182055644   \\
    \mathbf{0.44721359}6792237 &\mathbf{0.447213}596447341
   \end{array}
 \right)\;,
\end{equation*}
\begin{equation*}
  A_\mathbf{P} =\left(
   \begin{array}{ccccc}
   0.894427190999916 & 0.894427182055644 \\
    \mathbf{0.447213595}582167 & \mathbf{0.44721361338850}2
   \end{array}
 \right)\;.
\end{equation*}
It can be seen that $A_\mathbf{P}$ is a better estimate as it recovered more digits of the ground truth.  The mixtures are plotted in Fig. \ref{eg1_mix}, and the results are presented in Fig. \ref{eg1_rec}.

\begin{figure}
\includegraphics[height=5cm,width=8cm]{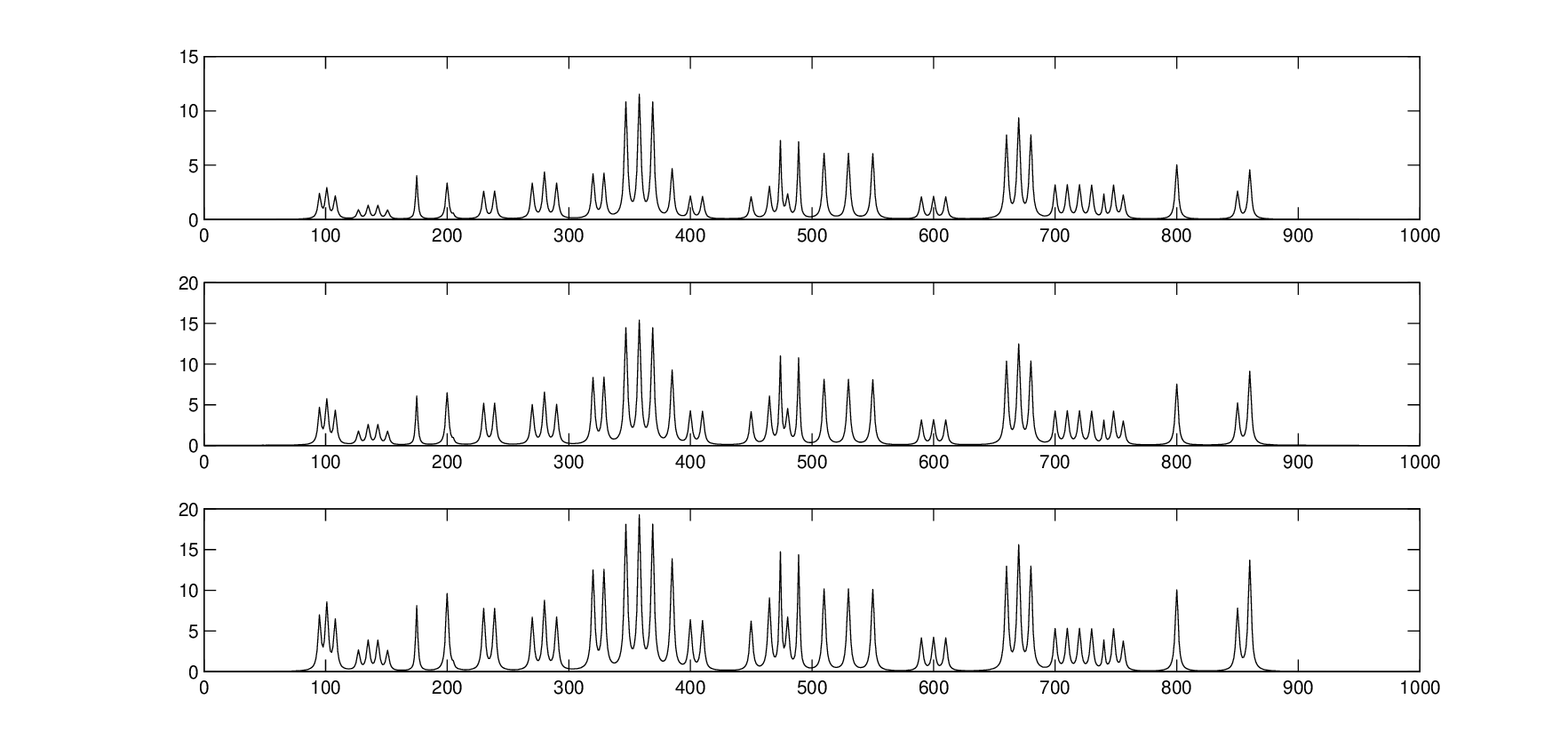}
\includegraphics[height=5cm,width=8cm]{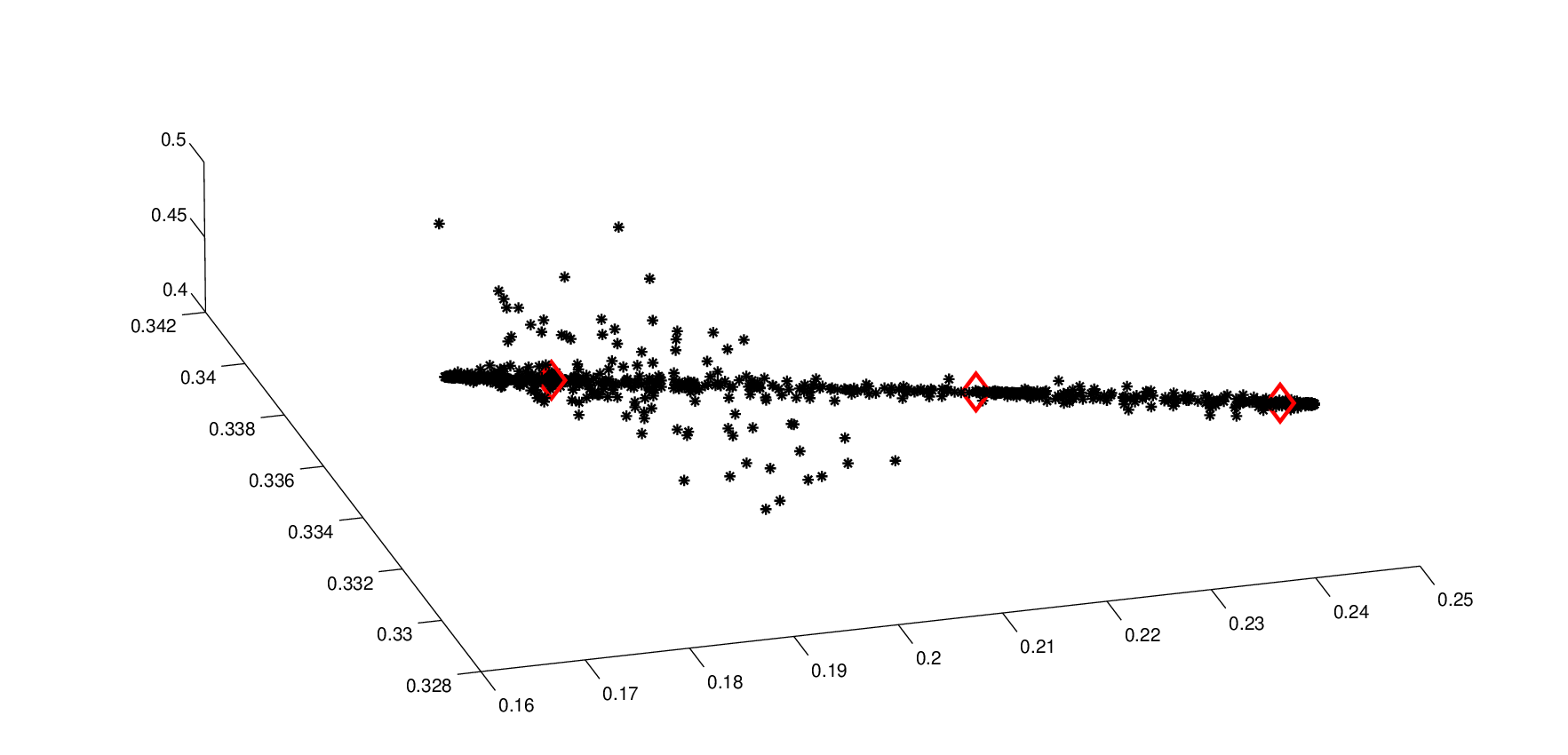}
\caption{three mixtures (left column) and their scattered plot (right column).}
\label{eg2_mix}
\end{figure}

\begin{figure}
\includegraphics[height=5cm,width=8cm]{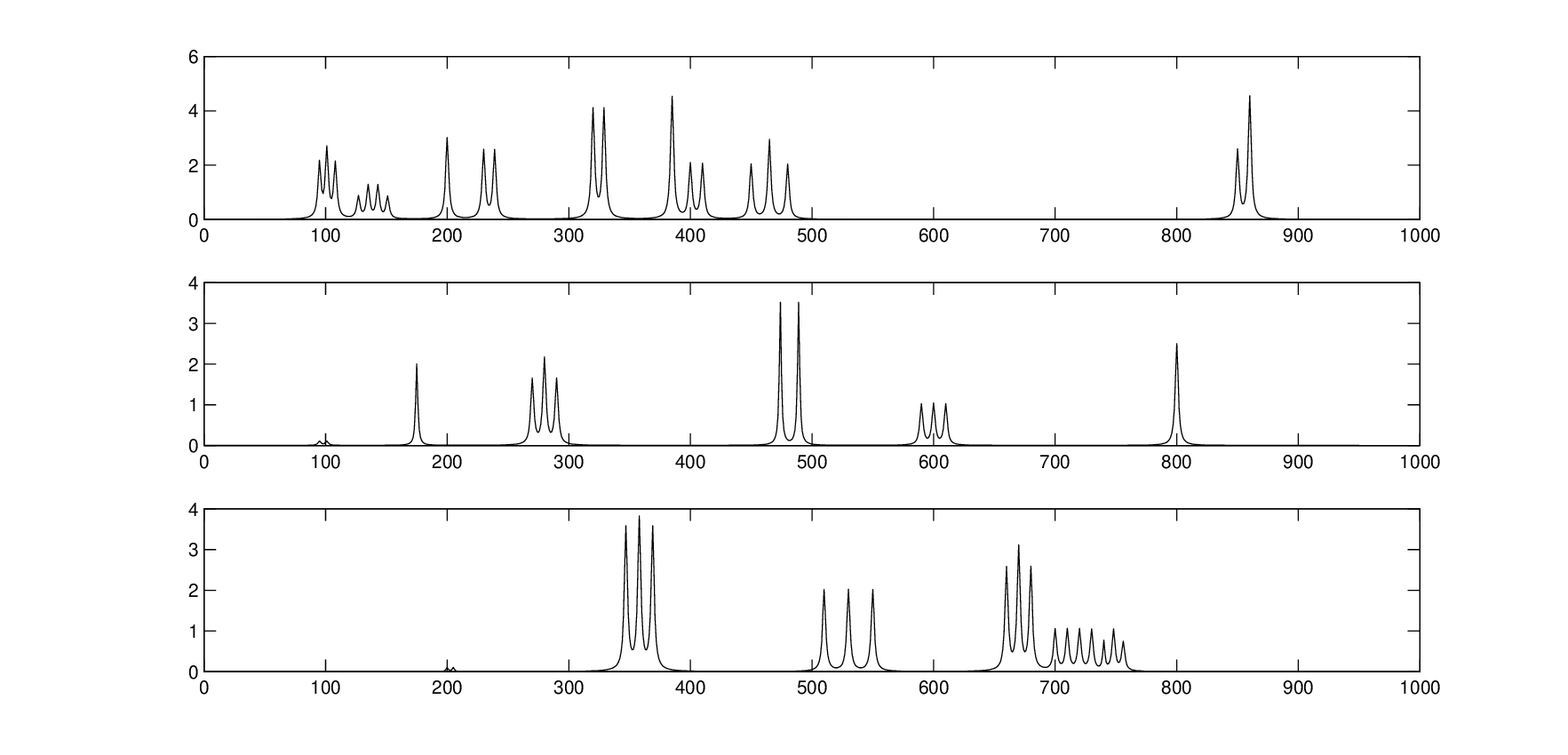}
\includegraphics[height=5cm,width=8cm]{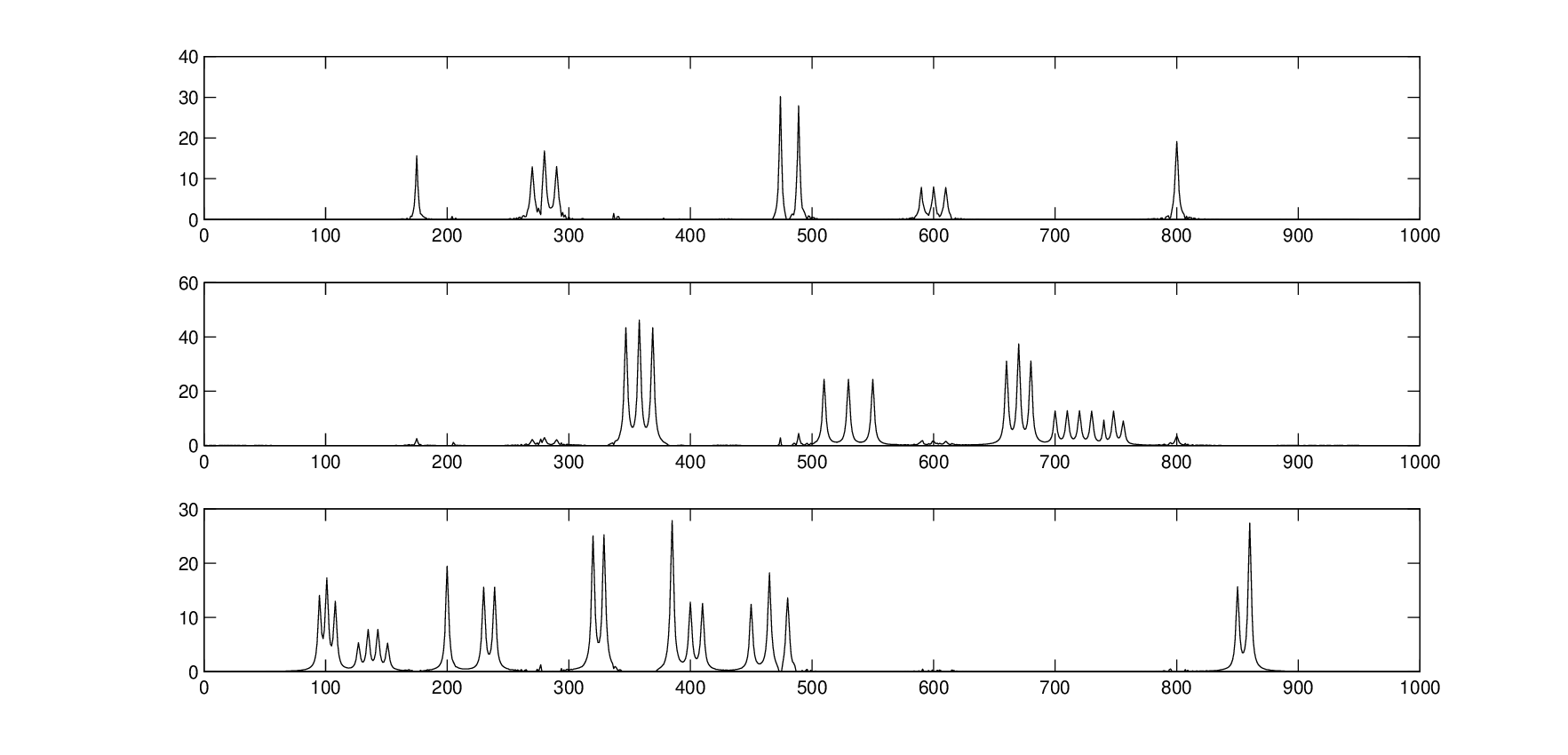}
\caption{three sources (left column) and their recovery by clustering and $\ell_1$ minimization (right column).}
\label{eg2_res}
\end{figure}

In the second example, three sources are to be separated from three mixtures.  The mixing matrix satisfies the OCDC condition, i.e., one of its columns is a nonnegative linear combination of the other two.  To test the robustness of the method, we added Gaussian noise (SNR = 60 dB) to the data.   The mixtures and their geometric structure are plotted in Fig. \ref{eg2_mix}. First the data clustering was used to obtain an estimation of the mixing matrix, then an $\ell_1$ optimization problem is solved to retrieve the sources. The results are shown in Fig. \ref{eg2_res}.  It can be seen that the recovered sources agree well with the ground truth.

\begin{figure}
\includegraphics[height=5cm,width=7.5cm]{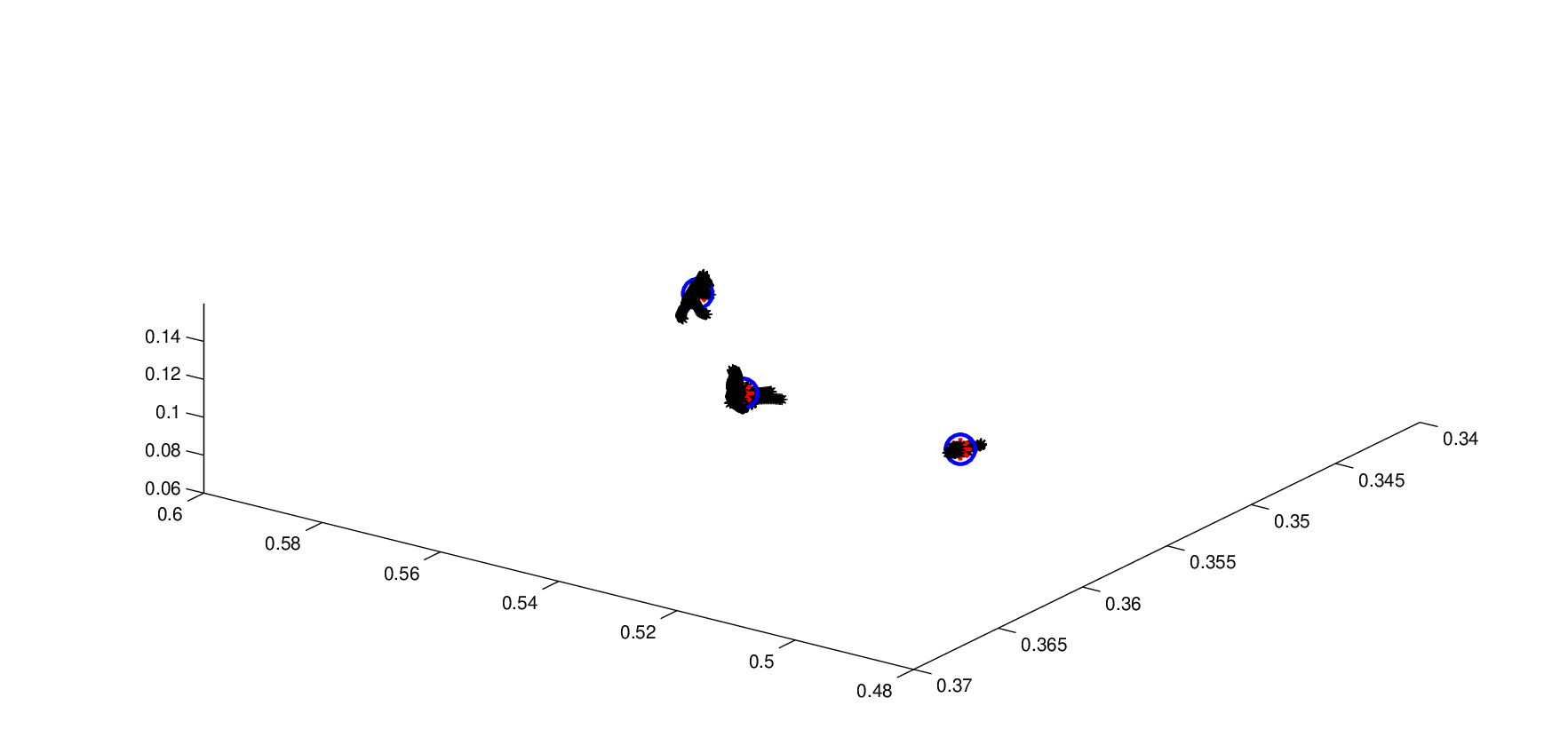}
\includegraphics[height=5cm,width=8cm]{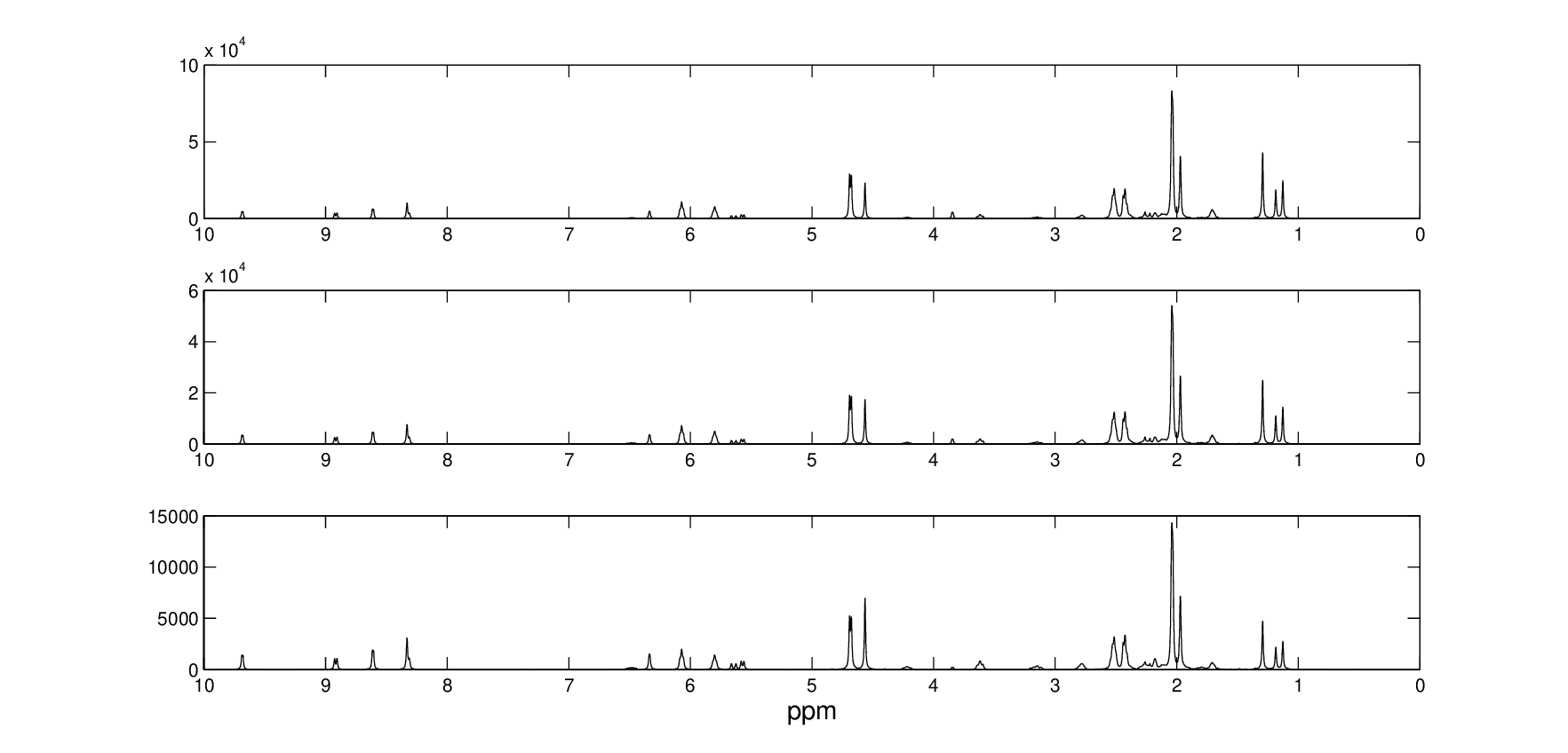}
\caption{three columns of $A$ are identified as the three center points in blue circles attracting most points in scatter plots of the columns of X (left), and the three rows of X (right).}
\label{real_data_3}
\end{figure}

\begin{figure}
\includegraphics[height=7cm,width=8cm]{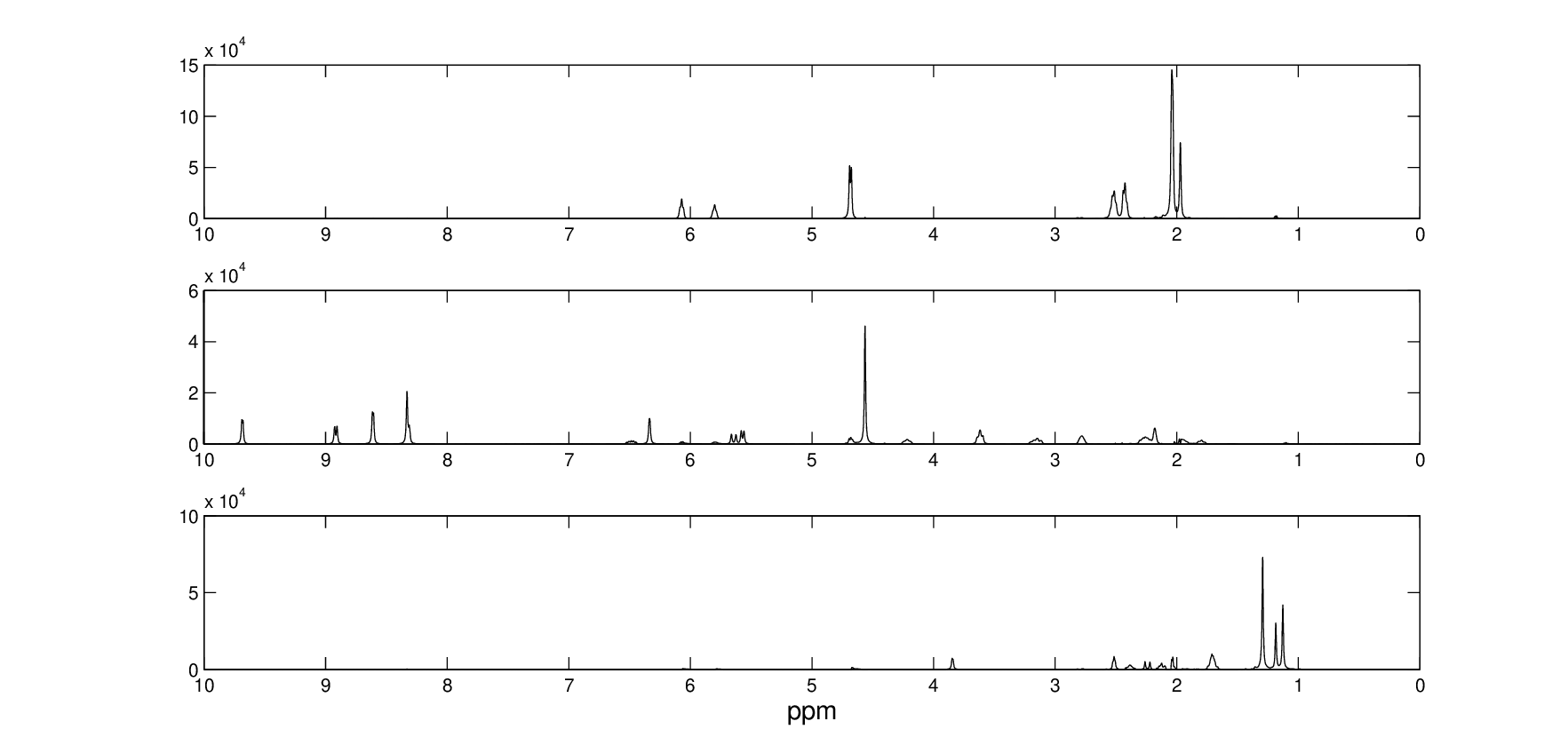}
\includegraphics[height=7cm,width=8cm]{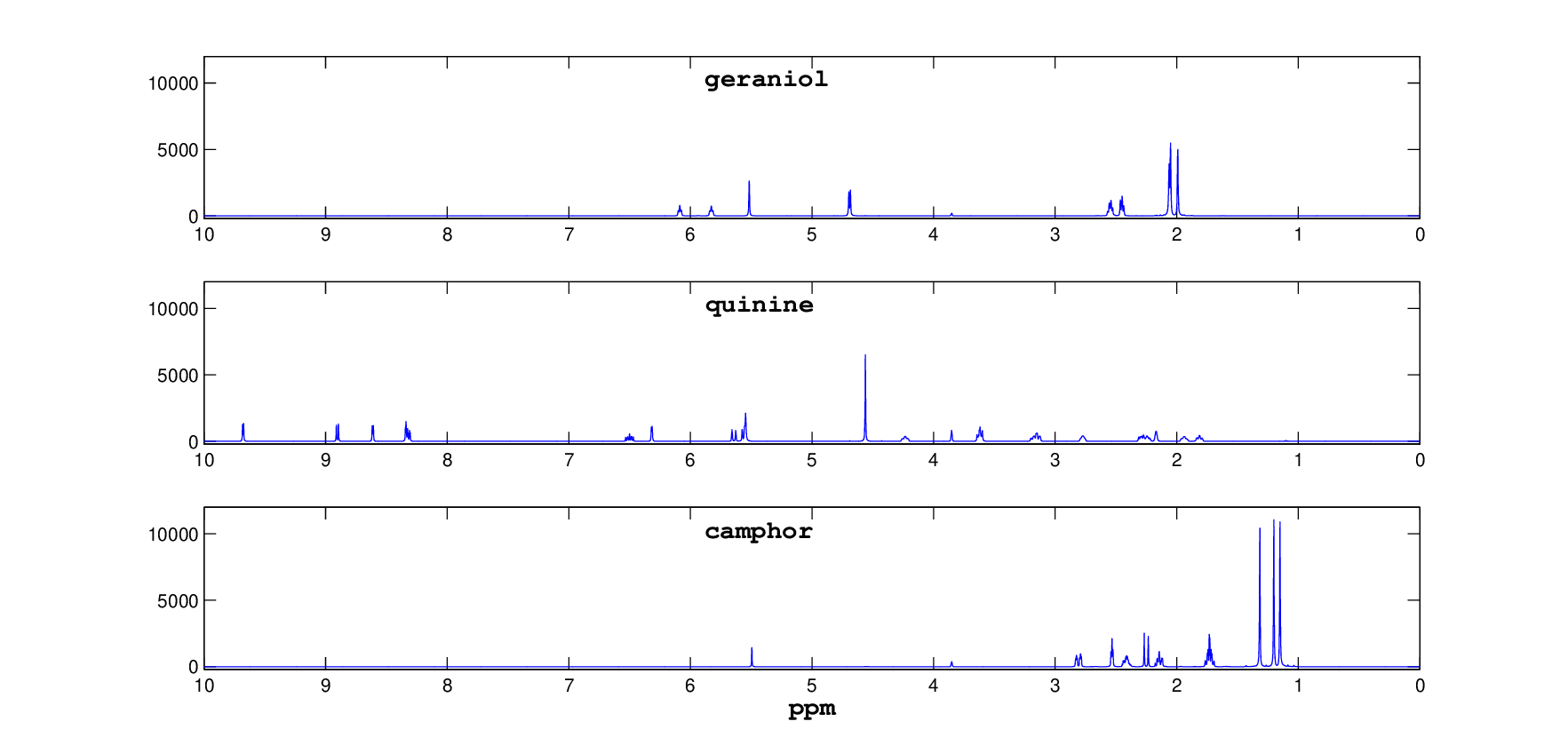}
\caption{the recovered source signals by nonnegative $\ell_1$ (left) and the ground truth (right).}
\label{real_result}
\end{figure}

\begin{figure}
\includegraphics[height=7cm,width=8cm]{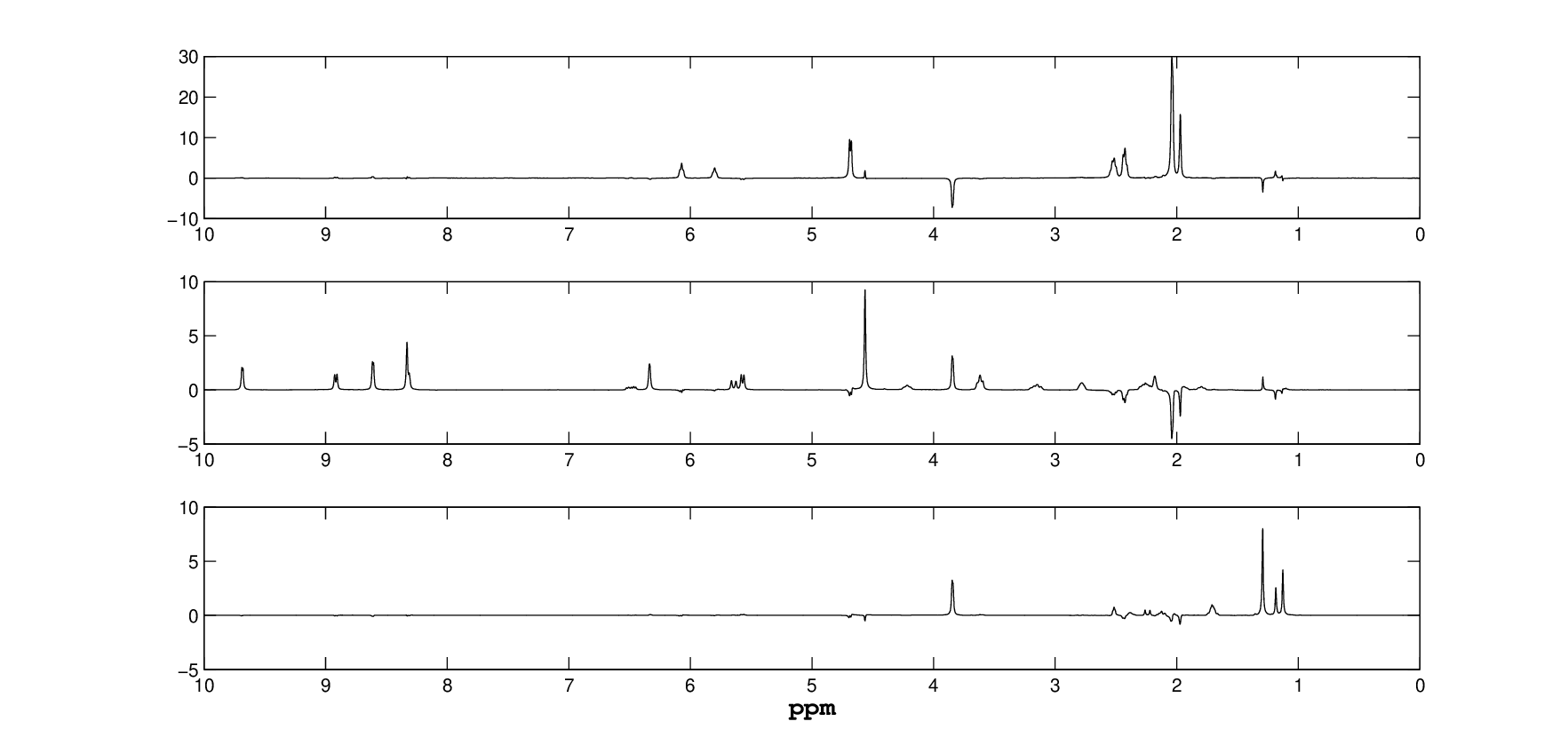}
\includegraphics[height=7cm,width=8cm]{figure/real_reference.eps}
\caption{the recovered source signals using NN method (left) and the ground truth (right).}
\label{real_result_NN}
\end{figure}

For the third example, we provide a set of real data to test our method.  The data is produced by diffusion ordered
spectroscopy (DOSY) which is an NMR spectroscopy technique used by chemists for
mixture separation \cite{Mor}.  However, the three compounds used in the experiment (quinine, geraniol, and
camphor) have similar chemical functional groups (i.e. there is overlapping in their NMR spectra) \cite{Nil},
for which DOSY fails to separate them.  It is known that each of the three sources has
dominant interval(s) over others in its NMR spectrum.   This can also be verified
from the three isolated clusters formed in their mixed NMR spectra (see the geometry of their
mixtures in Fig. \ref{real_data_3}).  Here we separate three sources from three mixtures.
Fig. \ref{real_data_3} plots the mixtures (rows of $X$) and their geometry (columns of $X$) where
three clusters of points can be spotted.  Then the columns of $A$ are identified as the
center points of three clusters.  The solutions are presented in Fig. \ref{real_result}, the results are satisfactory comparing with the ground truth.  As a comparison, the source signals recovered
by NN \cite{NN05} is shown in Fig. \ref{real_result_NN} where
$S = \mathrm{inverse}(A)\,X$, here the inverse is Moore-Penrose (the least squares sense) pseudo-inverse
which produces some negative (erroneous) peaks in $S$.


\subsection{Structure Assisted NMF examples}
We proposed nonnegative matrix factorization with two different cost functions to address the nearly degeneracy in the mixing matrix.  The first one imposed a penalty term on the differences of its columns (termed CD-NMF); while the second cost function addresses the degeneracy of the mixing matrix by adding a penalty on the inner products of all its columns (termed CP-NMF).  The comparison of the recovered mixing matrix by CD-NMF method,  CP-NMF method, and the ground truth are shown here
\begin{equation*}
  A_{\mathrm{CD-NMF}} =\left(
   \begin{array}{cccc}
   0.230364553564128 &  0.222325645752675 &  0.197973272413603\\
   0.428178141700514 &  0.426117653427851 &  0.447928422731077 \\
   0.709793002863661  & 0.677723684276001 &  0.563864889947392
   \end{array}
 \right)\;,\;
 \end{equation*}
 \begin{equation*}
A_{\mathrm{CP-NMF}} =\left(
   \begin{array}{cccc}
    0.242865565572616 &   0.180098969354018 &  0.285524082512094\\
   0.516635045464442  &   0.309883686246381 &   0.582543757260167\\
   0.729813101379756  &   0.538325711121780 &  0.857017001731787
      \end{array}
    \right)\;,\;
     \end{equation*}
\begin{equation*}
     A_{\mathrm{TR}} =\left(
   \begin{array}{cccc}
   1.000013203688969  & 1.099999061455011 &  3.999999720965063\\
   2.000007131698644  & 2.200024874750488 &  8.000009370903321 \\
   3.000004644386780  & 3.299985225130326 &  12.000000027664264

 \end{array}
    \right)\;
 \end{equation*}
The condition number of the ground truth matrix $A_{\mathrm{TR}}$ is about $5.56\times 10^6$, we actually started with a singular matrix
$
    \left(
   \begin{array}{cccc}
    1 &  1.1 &  4\\
    2  &  2.2 &  8 \\
    3  &  3.3 &  12
 \end{array}
    \right)\;,
$ then added Gaussian noise of SNR = 100. It is possible to compare the NMF results with the ground truth by scaling the rows of the matrices as in the first example.  Here we shall propose a metric to measure the distance of the matrices. The following distance between column degenerate matrices is proposed \footnote{The Comon's index\cite{Comon} used for measure the distance of nonsingular matrices in BSS problems is not suitable metric for nearly degenerate matrices.}
\begin{defi}
Consider two nearly degenerate matrices (parallel columns) $A$ and $\bar{A}$ of size $(m,n)$ with normalized columns.  The distance between $A$ and $\hat{A}$ denoted by $\Delta(A,\bar{A})$ is
\begin{equation*}
\Delta(A,\bar{A}) =\|A^{\mathbf{T}}\bar{A} -J_n\|_\mathbf{F}\;,
\end{equation*} where $J_n$ is a square matrix of size $n$ of ones.
\end{defi}Apparently the smaller distance is, the more similar of the matrices.  We computed the following matrix distances.
\begin{equation*}
\Delta(A_{\mathrm{TR}},A_{\mathrm{CD-NMF}}) = 0.005850257335153\;,\;  \Delta(A_{\mathrm{TR}},A_{\mathrm{CP-NMF}}) = 0.003635218995144.
\end{equation*} and the regularization parameters $\alpha = \beta = 0.001$.

\section{Conclusion}
This paper presented novel methods to retrieve source signals from the nearly degenerate mixtures.
The motivation comes from NMR spectroscopy of chemical compounds with similar diffusion rates.  Inspired by the NMR structure of these chemicals, we propose a viable source condition which requires dominant interval(s) from each source signal over the others.  This condition is well suited for many real-life signals.   Besides, the nearly degenerate mixtures are assumed to be generated from the following mixing process: either all the columns of the mixing matrix are parallel or one column is a nonnegative linear combinations of others.   We first use data clustering to identify the mixing matrix, then we develop two approaches to improve source signals' recovery.  The first approach minimizes a constrained quadratic program for a better mixing matrix, while the second method seeks the sparsest solution for each column of the source matrix by solving an $\ell_1$ optimization.  If no(or very limited) information on the source signals are available, two NMF variants are proposed by adding regularization terms to enforce the degeneracy of the columns, hence a desired solution can be obtained.  Numerical results on NMR spectra data show satisfactory performance of our method and offer promise towards understanding and detecting complex chemical spectra.  Though the methods are motivated by the NMR spectroscopy, the underlying ideas may be generalized to different data sets in other applications. For future work, we plan to investigate a mixture data separation problem where the mixing matrix has the following form
\[
\begin{bmatrix}
\exp(-\lambda D_1g_1^2) & \exp(-\lambda D_2g_1^2) & \dots& \exp(-\lambda D_ng_1^2)    \\
 \exp(-\lambda D_1g_2^2)&\exp(-\lambda D_2g_2^2)   & \dots& \exp(-\lambda D_ng_2^2)   \\
\vdots                   &\vdots                    & \vdots &  \vdots  \\
\exp(-\lambda D_1g_m^2) & \exp(-\lambda D_2g_m^2) & \dots& \exp(-\lambda D_ng_m^2)    \\
\end{bmatrix}
\]
Here $g$ is a controllable parameter, and $g_i's$ are linearly dependent.
Hence the entries of each column have nonlinear relations. We plan to study
how to impose constraints to guide the NMF solution revealing the nonlinearity among the entries.
\section*{Acknowledgements}  The authors wish to thank Professor A.J. Shaka and his group for their experimental NMR data. We also wish to thank Dr. Bin Yuan for helpful discussions.  YS was partially supported by Simons Foundation Grant 800006981.
JX was partially supported by NSF grant DMS-1924548.

\bibliographystyle{dcu}
\bibliography{Sun_NMF}

@article {Boardman_93,
    AUTHOR = {Boardman,J},
     TITLE = {Automated spectral unmixing of AVRIS data using convex geometry concepts},
   JOURNAL = {},
    VOLUME = {1},
      YEAR = {1993},
    NUMBER = {suppl. 1},
     PAGES = {11-14},
}

@article {ManiNMF,
    AUTHOR = {Cai,D and He, X.F and Han, J.W. and Huang, T.S},
     TITLE = {Graph Regularized Nonnegative Matrix Factorization for Data Representation},
   JOURNAL = {IEEE Transactions on Pattern Analysis and Machine Intelligence},
    VOLUME = {33},
      YEAR = {2011},
    NUMBER = {8},
     PAGES = {1548-1560},
}

@article {CanT,
    AUTHOR = {Cand\'{e}s, E. and  Romberg, J. and Tao,T.},
     TITLE = {Robust uncertainty principles: exact signal reconstruction from
              highly incomplete frequency information},
   JOURNAL = {IEEE Trans. Inform. Theory},
    VOLUME = {52},
      YEAR = {2006},
    NUMBER = {},
     PAGES = {489-509},
}

@book {Chang_07,
    AUTHOR = {Chang, C-I.},
     TITLE = {Hyperspectral Data Exploitation: Theory and Applications},
    SERIES = {},
    VOLUME = {},
 PUBLISHER = {Wiley-Interscience},
   ADDRESS = {},
      YEAR = {2007},
     PAGES = {},
      ISBN = {},
}

@article {Choi,
    AUTHOR = {Choi,S and Cichocki,A and Park,H and Lee,S},
     TITLE = {Blind source separation and independent component analysis: A review},
   JOURNAL = {Neural Inform. Process. Lett. Rev.},
    VOLUME = {6},
      YEAR = {2005},
    NUMBER = {},
     PAGES = {1-57},
}

@Book {Cic,
    AUTHOR = {Cichocki,A. and Amari,S.},
     TITLE = {Adaptive Blind Signal and Image Processing:
              Learning Algorithms and Applications},
    SERIES = {},
    VOLUME = {},
 PUBLISHER = {John Wiley and Sons},
   ADDRESS = {New York},
      YEAR = {2005},
     PAGES = {},
      ISBN = {},
}

@article {Comon,
    AUTHOR = {Comon,P.},
     TITLE = {Independent component analysis--a new concept ?},
   JOURNAL = {Signal Processing},
    VOLUME = {36},
      YEAR = {1994},
    NUMBER = {},
     PAGES = {287-314},
}

@Book {Comon1,
    AUTHOR = {Comon,P. and Jutten,C.},
     TITLE = {Handbook of Blind Source Separation: Independent Component
              Analysis and Applications},
    SERIES = {},
    VOLUME = {},
 PUBLISHER = {Academic Press},
   ADDRESS = {},
      YEAR = {2010},
     PAGES = {},
      ISBN = {},
}

@article {MVT,
    AUTHOR = {Craig,M.},
     TITLE = { Minimum-volume transformation for remotely sensed data},
   JOURNAL = {IEEE Transcations on Geoscience and Remote Sensing},
    VOLUME = {32},
      YEAR = {1994},
    NUMBER = {},
     PAGES = {542-552},
}

@article {Don,
    AUTHOR = {Donoho,D. and Tanner,J.},
     TITLE = {Sparse nonnegative solutions of underdetermined
              linear equations by linear programming},
   JOURNAL = {Proc Natl Acad Sci USA},
    VOLUME = {102},
      YEAR = {2005},
    NUMBER = {},
     PAGES = {9446-9451},
}

@article {G_O_09,
    AUTHOR = {Guo,Z. and Osher,S.},
     TITLE = { Template matching via $\ell_1$ minimization and its application
               to hyperspectral target detection},
   JOURNAL = {Tech. Rep., UCLA},
    VOLUME = {103},
      YEAR = {2009},
    NUMBER = {},
     PAGES = {},
}

@article {Lee,
    AUTHOR = {Lee,D.D. and Seung,H.S.},
     TITLE = {Learning of the parts of objects by non-negative matrix
              factorization},
   JOURNAL = {Nature},
    VOLUME = {401},
      YEAR = {1999},
    NUMBER = {},
     PAGES = {788-791},
}

@article {MVSA,
    AUTHOR = {Li,J. and Bioucas-Dias,J.M.},
     TITLE = { Minimum volume simplex analysis: a fast algorithm to unmix hyperspectral data},
   JOURNAL = {Geoscience and Remote Sensing Symposium},
    VOLUME = {3},
      YEAR = {2008},
    NUMBER = {},
     PAGES = {III250 - III 253},
}

@Book {Mor,
    AUTHOR = {Morris,G.},
     TITLE = {},
    SERIES = {},
    VOLUME = {},
 PUBLISHER = {John Wiley},
   ADDRESS = {New York},
      YEAR = {2001},
     PAGES = {},
      ISBN = {},
}

@article {NN05,
    AUTHOR = {Naanaa,W. and Nuzillard,J.M.},
     TITLE = {Blind source separation of positive and partially correlated data},
   JOURNAL = {Signal Processing},
    VOLUME = {85},
      YEAR = {2005},
    NUMBER = {9},
     PAGES = {1711-1722},
}

@article {VCA,
    AUTHOR = {Nascimento,J.M.P. and Bioucas-Diasm,J.M.},
     TITLE = { Vertex component analysis: a fast algorithm to unmix hyperspectral data},
   JOURNAL = {IEEE Transactions on Geoscience and Remote Sensing},
    VOLUME = {43},
      YEAR = {2005},
    NUMBER = {4},
     PAGES = {898-910},
}

@article {Nil,
    AUTHOR = {Nilsson,M. and Connel,M. and Davies,A. and Morris,G.},
     TITLE = {Biexponential Fitting of Diffusion-Ordered NMR
              Data: Practicalities and Limitations},
   JOURNAL = {Analytical Chemistry},
    VOLUME = {78},
      YEAR = {2006},
    NUMBER = {},
     PAGES = {3040-3045},
}

@article {NMF0,
    AUTHOR = {Paatero,P. and Tapper,U.},
     TITLE = {Positive matrix factorization: A non-negative factor
    model with optimal utilization of error estimates of data values},
   JOURNAL = {Environmetr.},
    VOLUME = {5},
      YEAR = {1994},
    NUMBER = {2},
     PAGES = {111-126},
}

@article {sun_xin_pNN,
    AUTHOR = {Sun,Y. and  Ridge,C. and del Rio, F. and Shaka,A.J. and Xin,J.},
     TITLE = {Postprocessing and sparse Blind source separation of
              positive and partially overlapped data},
   JOURNAL = {Signal Processing},
    VOLUME = {91},
      YEAR = {2011},
    NUMBER = {8},
     PAGES = {1838-1851},
}

@article {Winter_99,
    AUTHOR = {Winter,M.E.},
     TITLE = {N-findr: an algorithm for fast autonomous spectral
              endmember determination in hyperspectral data},
   JOURNAL = {Proc. of the SPIE},
    VOLUME = {3753},
      YEAR = {1999},
    NUMBER = {},
     PAGES = {266-275},
}

@article {YO,
    AUTHOR = {Yin,W. and Osher,S. and Goldfarb,D. and Darbon,J.},
     TITLE = {Bregman iterative algorithm for $\ell_1$-minimization with
              applications to compressive sensing},
   JOURNAL = {SIAM J. Image Sci},
    VOLUME = {1},
      YEAR = {2008},
    NUMBER = {},
     PAGES = {143-168},
}

@article {Zhang,
    AUTHOR = {Zhang,Y.},
     TITLE = {Theory of compressive sensing via $L_1$-Minimization:
              A Non-RIP analysis and extensions},
   JOURNAL = {Technical report, Rice University},
    VOLUME = {},
      YEAR = {2009},
    NUMBER = {},
     PAGES = {},
}

\end{document}